%% file: main.tex
\documentclass[10pt]{article}
\pdfoutput=1 
\usepackage[top=1.3in,bottom=1.2in,left=1in,right=1in]{geometry}
\usepackage{amsmath}
\usepackage{amssymb}
\usepackage{amsthm}
\usepackage{mathtools}
\usepackage[only,llbracket,rrbracket]{stmaryrd}
\usepackage[english]{babel}
\usepackage{graphicx}
\usepackage[font=small,figurewithin=section]{caption}
\usepackage{subcaption}
\usepackage{float}
\usepackage[numbers]{natbib} 
\usepackage{hyperref}
\usepackage[all]{xy}
\usepackage{multirow}
\usepackage{array}
\usepackage{enumitem}
\usepackage{tikz}
\usetikzlibrary{matrix,arrows,decorations.pathmorphing}
\usepackage{color}
\usepackage{authblk}
\usepackage[mathscr]{euscript}

\numberwithin{equation}{section}
\makeatletter
\g@addto@macro\bfseries{\boldmath}
\makeatother

\hypersetup{colorlinks=true,citecolor=black,filecolor=black,linkcolor=black,urlcolor=black}
\hypersetup{pdfstartview=FitB,pdfpagemode=UseNone}

\setlist{nolistsep}

\newcolumntype{L}[1]{>{\raggedright\let\newline\\\arraybackslash\hspace{0pt}}m{#1}}
\newcolumntype{C}[1]{>{\centering\let\newline\\\arraybackslash\hspace{0pt}}m{#1}}
\newcolumntype{R}[1]{>{\raggedleft\let\newline\\\arraybackslash\hspace{0pt}}m{#1}}
\newcolumntype{N}{@{}m{0pt}@{}}

\SetSymbolFont{stmry}{bold}{U}{stmry}{m}{n}

\input{preamble}

\begin{document}


\title{\vspace*{-35pt}High-order polygonal discontinuous Petrov-Galerkin (PolyDPG) methods using ultraweak formulations}
\author[1,2]{Ali Vaziri Astaneh\thanks{Corresponding author. E-mail: ali@vaziri.info}}
\author[1]{Federico Fuentes}
\author[1]{Jaime Mora}
\author[1]{Leszek Demkowicz}
\affil[1]{The Institute for Computational Engineering and Sciences (ICES), The University of Texas at Austin, 201 E 24th St, Austin, TX 78712, USA}
\affil[2]{MSC Software Corporation, Newport Beach, CA 92660, USA}
\date{\vspace{-10mm}}

\maketitle
\input{abstract}

\input{1_Introduction}
\input{2_Method}

\input{3_Results}
\input{4_Conclusions}

\appendix
\input{5_Convergence}

\bibliographystyle{apalike}
\bibliography{main}


\end{document}

%% file: preamble.tex
\newcommand{\R}{\mathbb{R}} 
\newcommand{\N}{\mathbb{N}} 
\newcommand{\HSo}{H} 
\newcommand{\bHSo}{\boldsymbol{H}} 
\newcommand{\Leb}{L} 
\newcommand{\bLeb}{\boldsymbol{L}} 

\newcommand{\idmatrix}{\textup{\uppercase\expandafter{\romannumeral 1}}}
\newcommand{\mesh}{\mathcal{T}} 
\DeclareMathOperator{\dd}{d} 
\newcommand{\T}{\mathsf{T}} 
\newcommand{\bdry}{\partial} 
\newcommand{\tr}{\mathrm{tr}} 
\newcommand{\opt}{\mathrm{opt}} 
\newcommand{\nopt}{\mathrm{n\text{-}opt}} 
\newcommand{\overlap}{\mathrm{ov}} 

\newcommand{\btau}{\boldsymbol{\tau}}

\DeclareMathOperator{\diam}{diam}
\DeclareMathOperator*{\argmin}{arg\,min}

\DeclareMathOperator{\curl}{curl}
\let\div\relax 
\DeclareMathOperator{\div}{div}

\newtheoremstyle{boldremark}
    {\dimexpr\topsep/2\relax} 
    {\dimexpr\topsep/2\relax} 
    {}          
    {}          
    {\bfseries} 
    {.}         
    {.5em}      
    {}          
\theoremstyle{definition}
\newtheorem*{definition*}{Definition}
\theoremstyle{plain}
\newtheorem{theorem}{Theorem}
\numberwithin{theorem}{section}

\numberwithin{lemma}{section}

\numberwithin{corollary}{section}
\theoremstyle{boldremark}
\newtheorem{remark}{Remark}
\numberwithin{remark}{section}

\makeatletter
\newcommand{\dashto}[1][2pt]{
  \settowidth{\@tempdima}{${}\rightarrow{}$}
  \makebox[\@tempdima]{${}\rightarrow{}$}
  \makebox[-\@tempdima]{\hspace{-0.1\@tempdima}\color{white}\rule[0.5ex]{#1}{1pt}}
  \makebox[\@tempdima]{}
  }
\makeatother
\let\tilde\widetilde

\newcommand{\upperone}[2]{\raisebox{-0.2\height}{$#1{}^1\!$}}
\newcommand{\upone}{\mathpalette\upperone\relax}
\newcommand{\lowertwo}[2]{\raisebox{-0.1\height}{$#1{}_{\!2}$}}
\newcommand{\lowtwo}{\mathpalette\lowertwo\relax}

\newcommand{\frachalf}[2]{\raisebox{0.0\height}{$#1{}\upone/\lowtwo$}}
\newcommand{\onehalf}{\mathpalette\frachalf\relax}

\newcommand{\mcE}{\mathcal{E}}

\newcommand{\mcK}{\mathcal{K}}

\newcommand{\mcO}{\mathcal{O}}
\newcommand{\mcP}{\mathcal{P}}
\newcommand{\mcQ}{\mathcal{Q}}
\newcommand{\mcR}{\mathcal{R}}

\newcommand{\mcT}{\mathcal{T}}

\newcommand{\sfu}{\mathsf{u}}

\newcommand{\sfB}{\mathsf{B}}

\newcommand{\sfG}{\mathsf{G}}

\newcommand{\bsq}{\boldsymbol{q}}

\newcommand{\scB}{\mathscr{B}}

\newcommand{\scP}{\mathscr{P}}

\newcommand{\scR}{\mathscr{R}}

\newcommand{\scT}{\mathscr{T}}
\newcommand{\scU}{\mathscr{U}}
\newcommand{\scV}{\mathscr{V}}
\newcommand{\scW}{\mathscr{W}}

\newcommand{\fkh}{\mathfrak{h}}

\newcommand{\fku}{\mathfrak{u}}
\newcommand{\fkv}{\mathfrak{v}}

\newcommand{\fkH}{\mathfrak{H}}

\newcommand{\overbar}[1]{\mkern 1.5mu\overline{\mkern-1.5mu#1\mkern-1.5mu}\mkern 1.5mu}

\newcommand\bovermat[2]{\makebox[0pt][l]{$\smash{\overbrace{\phantom{
	\begin{matrix}#2\end{matrix}}}^{#1}}$}#2}

\newcommand{\nml}{\hat{\mathbf{n}}}
\newcommand{\sfellt}[2]{\raisebox{0.0\height}{$#1{}\text{\textsf{\fontfamily{PTSans-TLF}\selectfont l}}$}}
\newcommand{\sfell}{\mathpalette\sfellt\relax}

%% file: abstract.tex
\renewcommand{\abstractname}{\large Abstract}
\begin{abstract}
\small 
\noindent

This work represents the first endeavor in using ultraweak formulations to implement high-order polygonal finite element methods via the discontinuous Petrov-Galerkin (DPG) methodology.
Ultraweak variational formulations are nonstandard in that all the weight of the derivatives lies in the test space, while most of the trial space can be chosen as copies of $\Leb^2$-discretizations that have no need to be continuous across adjacent elements.
Additionally, the test spaces are broken along the mesh interfaces.
This allows one to construct conforming polygonal finite element methods, termed here as PolyDPG methods, by defining most spaces by restriction of a bounding triangle or box to the polygonal element.
The only variables that require nontrivial compatibility across elements are the so-called interface or skeleton variables, which can be defined directly on the element boundaries.
Unlike other high-order polygonal methods, PolyDPG methods do not require ad hoc stabilization terms thanks to the crafted stability of the DPG methodology.
A proof of convergence of the form $h^p$ is provided and corroborated through several illustrative numerical examples.
These include polygonal meshes with $n$-sided convex elements and with highly distorted concave elements, as well as the modeling of discontinuous material properties along an arbitrary interface that cuts a uniform grid. 
Since PolyDPG methods have a natural \textit{a posteriori} error estimator a polygonal adaptive strategy is developed and compared to standard adaptivity schemes based on constrained hanging nodes. 
This work is also accompanied by an open-source \texttt{PolyDPG} software supporting polygonal and conventional elements.

\end{abstract}

%% file: 1_Introduction.tex
\section{Introduction}
\label{sec:Introduction}

Numerical solutions of boundary value problems with meshes of general polytopes were first proposed by Wachspress \cite{wachspress}, who introduced rational barycentric coordinates that formed a finite element basis over convex polygons, leading to a conforming finite element method (FEM) with new types of elements.
Over the last two decades, there has been a growing collection of numerical methods using general polytopes which extend well beyond the original ideas of Wachspress.
Among the reasons for this group of methods to thrive is a handful of advantages that polytopes offer over traditionally shaped elements (simplices, hexahedra, etc.).
These include: matching complex interfaces (see e.g.~\cite{wg1,vem_interface}); greater flexibility to mesh complex geometries and their role as transition elements \cite{pfem_conforming}; avoiding the limitations of parametric elements for highly distorted or ill-shaped elements (see e.g.~\cite{vem_finite_elast,parametricdistortion}); handling multiple hanging nodes in local $h$-refinements \cite{pfem_quadtree}; and allowing for greater deformations and less tendency to mesh-locking in incompressible media \cite{pfem_finite_elast}.


The features just mentioned give polytopal FEMs a wide range of applicability, especially where conventional methods do not fare well.
In fact, they are useful for resolving problems involving the deformation of materials with heterogeneous microstructure \cite{microstructure}, modeling complex materials like elastomers and biomaterials \cite{pfem_finite_elast,nem_bone}, creating meshes where interface fitting is required \cite{vem_interface}, and modeling fractured media \cite{vem_fracture_network}.  
Promising results have also been obtained in crack propagation modeling \cite{pfem_cracks,fracturemeshbias,Bishop2009,Bishop2016} and in topology optimization \cite{pfem_optimization,vem_optimization,antonietti2017,polytopmatlab}, since polygonal meshes combine the ability to mesh complex geometries with a reasonable number of elements while reducing mesh-induced bias in particular directions (which occurs in structured meshes of triangles or quadrilaterals) \cite{pfem_optimization,fracturemeshbias,antonietti2017}.

Many methods still utilize different types of generalized barycentric coordinates (including some valid in nonconvex polytopes), which have proliferated since Wachspress originally introduced them, as well as other choices of shape functions (see e.g.~\cite{Bishop2014}).
These methods are usually $\HSo^1$-conforming Galerkin FEMs \cite{pfem_conforming}, but there are some extensions to mixed methods (see e.g.~\cite{pfem_finite_elast}).
They mostly allow very flexible refinement schemes while avoiding constrained approximations \cite{pfem_quadtree}, but they are typically limited by first order $h$-convergence. 
Some families of high-order shape functions have been proposed, but only for convex polytopes (see e.g.~\cite{pfem_quadratic,pfem_families}). 
As the barycentric coordinates are in general rational polynomials, another challenge is the choice of the quadrature scheme used for integration \cite{pfem_integration1,pfem_integration2}.
Mimetic finite difference (MFD) methods are based on another discretization technique which also supports polygonal elements.
The technique consists of designing discrete differential operators such that fundamental vector calculus identities and physical laws can be reproduced in a discrete context \cite{mfd1,mfd2,mfd3}. 
Later, the ideas of MFDs led to the development of virtual element methods (VEMs) \cite{vem_basic}. 
In VEMs, appropriate spaces are tailored for each polytopal element, such that their functions have continuous and piecewise polynomial traces over the boundaries. The integrals over the cells can be computed exactly (i.e.~up to machine precision) with quadrature points only on the boundary \cite{pfem_vem_review}. 
The power of VEMs lies partly in eliminating the need of explicitly constructing the shape functions in the element, and yet resulting in a FEM-like variational setting \cite{vem_elliptic}. 
They are also high-order methods \cite{BeiraoHighOrderVEM}, and recent work has resulted in the construction of $\bHSo(\div)$- and $\bHSo(\curl)$-conforming spaces \cite{vem_hcurl_hdiv}. 
VEMs have been used for different problems like linear elasticity, plate bending, and second-order elliptic problems \cite{vem_elast,vem_plates,vem_elliptic}. 
But it must be noted that VEMs need a problem-dependent stability operator to guarantee their convergence \cite{pfem_vem_review}, and the solution at interior points of the elements is not accessible directly, so it has to be approximated \cite{vem_elliptic}.

Another method is the polytopal interior penalty $hp$ discontinuous Galerkin (IPDG) method \cite{hpdgfem}.
It is a nonconforming high-order method, which uses restrictions of standard FE spaces associated to a bounding box of each element. 
Due to its nonconformity, the method has a thorough but nonstandard equation-dependent error analysis, and like VEMs, it needs adding extra terms to ensure stability. 
Lastly, other recent methods include hybrid mimetic mixed methods \cite{gdmoriginal,hhmoriginal}, PFEM-VEM \cite{pfem_vem_review}, the weak Galerkin (WG) method \cite{wg1,wg2,wgoriginal}, hybrid high-order (HHO) methods \cite{HHOreview}, and hybridizable discontinuous Galerkin (HDG) methods \cite{HDGReview,HDG2D}.
More details on the historical development can be found in the thorough review \cite{pfem_vem_review}.

The objective of this article is to present a completely new family of high-order methods termed polygonal discontinuous Petrov-Galerkin (PolyDPG) methods.
They are based on so-called ``broken'' ultraweak variational formulations discretized using the discontinuous Petrov-Galerkin (DPG) methodology \cite{DPGOverview}.
These formulations, despite being well-defined at the infinite-dimensional level, admit a very large degree of discontinuities in both the trial and test spaces, since their test spaces are broken (i.e.~they may be discontinuous across element interfaces) and part of their trial spaces is in $\Leb^2$.
In fact, the only communication between elements happens through the so-called skeleton (or interface) variables that live on the element boundaries.
These nonstandard formulations can be systematically discretized in a conforming fashion (i.e., with discrete trial and test spaces that are subspaces of the infinite-dimensional ones) and solved using the variationally versatile DPG methodology, which always produces a positive-definite finite element stiffness matrix. 
The DPG methodology is essentially crafted to produce stability by using optimal test functions and without resorting to additional stabilization terms.
DPG methods have been successfully used for equations involving numerical stability issues \cite{cohen2012adaptivity,demkowicz2013robust,chan2014robust,Niemi2013,DLSFEM}, and applied to various physical problems such as wave propagation \cite{zitelli2011class,GopMugaHelmholtz,WavenumberHelmholtzDPG,petrides2016ices}, transmission problems \cite{HeuerKarkulikTransmission,HeuerFuhrerCoupling}, electromagnetism \cite{BrokenForms15}, elasticity \cite{Keith16,Bramwell12,FuentesCoupledElasticity,FuentesViscoelasticity}, fluid flow \cite{roberts2015discontinuous,ChanCompressible,ellisfluids,KeithOldroydB} and optical fibers via Schr\"{o}dinger's equation \cite{NagarajSchroedinger}.

In this paper we consider 2D problems, where the element boundaries are merely line segments, so high-order discretization of the skeleton variables is straightforward.
As we will show, this makes the broken ultraweak formulations an ideal framework for defining polygonal elements, and it results in the conforming FEMs we refer to as PolyDPG methods. 
PolyDPG methods are competitive with other existing polygonal methods, since they arise from very different ideas and they inherit many advantages from the DPG methodology.
For example, they can be easily generalized to different linear equations; they have a solid mathematical background in terms of proving stability and high-order convergence; they allow for discontinuous material properties while retaining stability; they result in positive-definite stiffness matrices; and they carry a completely natural arbitrary-order \textit{a posteriori} error estimator, which facilitates implementation of adaptive refinement strategies.
The last feature is particularly desirable when combined with polygonal elements, because there is no need for the constrained approximation technology to treat hanging nodes, paving the way for use in applications like dynamic fracture \cite{pfem_cracks,fracturemeshbias,Bishop2009,Bishop2016} and topology optimization \cite{pfem_optimization,vem_optimization,antonietti2017,polytopmatlab}.
We complement this article by providing an open-source software in MATLAB\textsuperscript{\textregistered}, also named \texttt{PolyDPG} \cite{polydpg}.

The outline of the article is as follows.
In Section \ref{sec:method} we describe a PolyDPG method for a model problem (Poisson's equation), along with the DPG solution scheme and the convergence theory (with the proof relegated to Appendix \ref{sec:AppConvergence}). 
In Section \ref{sec:results} several illustrative examples are presented.
High-order convergence for different $p$ is verified for both convex and highly distorted concave elements.
Then, a physically relevant problem involving discontinuous material properties along an arbitrary interface is solved.
Finally, an adaptive refinement strategy is described, successfully implemented, and compared to traditional adaptive schemes.
Our concluding remarks are presented in Section \ref{sec:conclusions}.

%% file: 2_Method.tex
\section{PolyDPG methods} 
\label{sec:method}

Typical FEMs map elements from the actual physical space to a known fixed master element space corresponding to the same element type.
For example, in 2D a general quadrilateral in $\R^2$ is mapped to a master quadrilateral (typically $(0,1)^2$ or $(-1,1)^2$).
This requires defining a master element for each element type, which is possible for limited types of elements (e.g. quadrilaterals and triangles in 2D, or hexahedra, tetrahedra, triangular prisms and pyramids in 3D), but is usually nonviable when dealing with general polytopes.
Thus, as with any polytopal FEM, the idea is to circumvent any master elements by shifting the focus directly to the physical space. 

The main issue in doing so is satisfying inter-element continuity of the basis functions, which is required for discretizing Sobolev spaces such as $H^1$.
This is partly resolved by using generalized barycentric coordinates, but these techniques are usually limited to first order methods (in terms of convergence), and it becomes difficult to discretize other Sobolev spaces such as $\bHSo(\curl)$ and $\bHSo(\div)$ even for the lowest order cases \cite{YanqiuPolytope}.
Indeed, even with the ``traditional'' pyramid element, having high-order discretizations for different spaces is challenging to achieve \cite{Nigam_Phillips_11,Fuentes2015,GillettePyramid,AinsworthPyramids}, and so is the case for 2D non-affine quadrilaterals \cite{arbogastquads}.
To overcome this, VEMs concentrate on the boundaries while nonconforming polytopal discontinuous methods, like IPDG, HHO, WG, and HDG (which are closely related \cite{BridgingHHOHDG,HDGReview}), remove the continuity requirements altogether.
However, all of these methods need to carefully add (equation-dependent) stabilization or penalty terms \cite{vem_basic,hpdgfem,HHOreview,wgoriginal,HDG2D}, and they must account for these in the error analysis, leading to a nonstandard theory of convergence \cite{Ciarlet94}.

As will be seen, the discontinuous Petrov-Galerkin (DPG) methodology is very general from a variational standpoint, so it is not limited to the traditional primal and mixed formulations.
Thus, without sacrificing any desirable stability properties, it is able to discretize ``broken'' ultraweak variational formulations, which avoid most inter-element continuity requirements.
The only continuity requirements are met by skeleton variables which live on the element boundaries.
Technically speaking, the resulting method is still a conforming FEM, and the ``standard'' error analysis can be applied.
This is very useful, because it allows to generalize the method to any well-posed linear equation formulated with traditional functional spaces ($\HSo^1$, $\bHSo(\curl)$, $\bHSo(\div)$ and $\Leb^2$).

In 2D, the polygonal element boundaries are simply line segments, so it is easy to define high-order discretizations along the mesh skeleton.
Given that this is less trivial for polyhedra in 3D, we only analyze 2D problems in this introductory paper.
We now proceed by introducing the model problem and its corresponding ultraweak formulations in the next section. 

\subsection{Model problem and ultraweak variational formulations} 

As a model problem, consider Poisson's equation coming from the steady-state heat equation in a (heterogeneous) domain $\Omega\subseteq\R^2$, where $u$ is the temperature, $\bsq$ is the heat flux, $k>0$ is the variable thermal conductivity, and $r$ is the internal heat source,
\begin{equation}
	-\div(k\nabla u)=r\,,\qquad\Leftrightarrow\qquad
	\left\{
  \begin{aligned}
  	\div\bsq &= r\,,\\
   	\textstyle{\frac{1}{k}}\bsq+\nabla u&= 0\,.
  \end{aligned}
 	\right.
	\label{eq:poissonheat}
\end{equation}
Note that the equation can be written directly as a second order system (left) or as a first order system (right).
For simplicity, we assume temperature boundary conditions along all of $\bdry\Omega$, so that $u=g$ at $\bdry\Omega$, where $g$ is a known function.

To solve the equation using FEMs, a variational form is required, and in this respect, there are many possibilities.
For now assume vanishing temperature boundary conditions so that $g=0$.
The classical approach stems directly from the second order equation by multiplying by a test function and integrating by parts once, leading to the primal formulation where the solution $u$ is sought in the trial space $\scU^\scP$ and must satisfy
\begin{equation}
	\begin{gathered}
		b^\scP(u,v)=\ell^\scP(v)\,\qquad\forall v\in \scV^\scP=\scU^\scP=\HSo_0^1(\Omega)\,,\\
		b^\scP(u,v)=(k\nabla u,\nabla v)_\Omega\,,\qquad\qquad \ell^\scP(v)=(r,v)_\Omega\,,
	\end{gathered}
	\label{eq:Primal}
\end{equation}
with $(u,v)_K=\int_K u\cdot v\dd K$ for $K\subseteq\Omega$.
Notice in this case $\scV^\scP=\scU^\scP$, so both spaces can be discretized in the same way, leading to the Galerkin method.
The same property holds for standard mixed formulations which stem from the first order system.
The ultraweak formulation is also derived from the first order system, but here all equations are integrated by parts to pass the derivatives to the test functions.
The resulting ultraweak formulation seeks $(u,\bsq)=\fku_0\in\scU_0=\Leb^2(\Omega)\times\bLeb^2(\Omega)$ satisfying
\begin{equation}
	\begin{gathered}
		b_0(\fku_0,\fkv_0)=\ell(\fkv_0)\,\qquad\forall(v,\btau)=\fkv_0\in\scV_0=\HSo_0^1(\Omega)\times\bHSo(\div,\Omega)\,,\\
		b_0\big((u,\bsq),(v,\btau)\big)=-(\bsq,\nabla v)_\Omega+(\textstyle{\frac{1}{k}}\bsq,\btau)_\Omega-(u,\div\btau)_\Omega,
			\qquad \ell\big((v,\btau)\big)=(r,v)_\Omega\,,
	\end{gathered}
	\label{eq:UW}
\end{equation}
where $\bLeb^2(\Omega)=(\Leb^2(\Omega))^2$.
Clearly the trial and test spaces in this case are completely different, $\scU_0\neq\scV_0$.
Thus, to solve this system it is necessary to drift away from the traditional Galerkin method.
As we will see, a discretization via minimum residual FEMs is a viable option.
It is worth remarking that the primal and ultraweak formulations are mutually well-posed in the infinite-dimensional setting \cite{Keith16,DemkowiczClosedRange,BrokenForms15}. 
Since the primal formulation is known to be well-posed in view of the Lax-Milgram theorem and Poincar\'{e}'s inequality, so is the ultraweak formulation.
This guarantees the existence of a unique solution in the trial space satisfying a stability estimate.

The ultraweak formulation has copies of $\Leb^2(\Omega)$ as a trial space, thus its discretization does not require satisfying any inter-element continuity, which is very desirable for polygons.
However, all the difficulties are passed to the test space for which inter-element continuity requirements are essential.
Fortunately, it is possible to remove these requirements in the test space as well, but at the cost of introducing skeleton variables, as we will see shortly. 
In fact, the practicality of DPG methods relies on using broken (or discontinuous) test spaces, and this results in a slightly modified formulation called the \textit{broken} ultraweak formulation, which will be derived in what follows.
Consider a mesh (i.e.~an open partition), $\mesh$, of $\Omega$ comprised of (disjoint) elements $K\in\mesh$, and define the broken spaces and piecewise integration,
\begin{equation}
	\begin{aligned}
		\HSo^1(\mesh)&=\{v\in\Leb^2(\Omega)\mid v|_K\in\HSo^1(K)\,,\forall K\in\mesh\}\,,\\
		\bHSo(\div,\mesh)&=\{\btau\in\bLeb^2(\Omega)\mid \btau|_K\in\bHSo(\div,K)\,,\forall K\in\mesh\}\,,\\
		{}&\qquad\quad(u,v)_\mesh=\sum_{K\in\mesh}(u|_K,v|_K)_K\,.
	\end{aligned}
\end{equation}
Then, \textit{element-wise}, multiply by broken test functions $(v,\btau)=\fkv\in\scV=\HSo^1(\mesh)\times\bHSo(\div,\mesh)$, integrate by parts, and sum across all elements.
The result is very similar to the ultraweak formulation, but has new terms on the boundaries of the elements involving $u|_{\bdry K}$ and $\bsq|_{\bdry K}\!\cdot\!\nml_K$, where $\nml_K$ is the outward normal to the element $K$.
These terms vanish if the test space is not broken (i.e.~$\scV_0$).
Unfortunately, if we want $u\in\Leb^2(\Omega)$ and $\bsq\in\bLeb^2(\Omega)$, then the traces $u|_{\bdry K}$ and $\bsq|_{\bdry K}\!\cdot\!\nml_K$ technically do not exist \cite{McLeanSobolev} and to incorporate them it is necessary to add new skeleton (or interface) variables in the spaces
\begin{equation}
	\begin{aligned}
		\HSo_0^{\onehalf}(\bdry\mesh)&=\big\{\hat{u}\mid\exists u\in\HSo_0^1(\Omega),\,\, 
			\hat{u}=\textstyle{\prod_{K\in\mesh}}(u|_K)\big|_{\bdry K}\big\}\,,\\
		\HSo^{-\onehalf}(\bdry\mesh)&=\big\{\hat{q}_{\nml}\mid\exists\bsq\in\bHSo(\div,\Omega),\,\, 
			\hat{q}_{\nml}=\textstyle{\prod_{K\in\mesh}}(\bsq|_K)\big|_{\bdry K}\!\cdot\!\nml_K\big\}\,,\\
		{}&\quad\langle \hat{u},\hat{v}\rangle_{\bdry\mesh}=\sum_{K\in\mesh}\langle(\hat{u})_K,(\hat{v})_K\rangle_{\bdry K}\,,
	\end{aligned}
	\label{eq:interfacespaces}
\end{equation}
where the duality $\langle\cdot,\cdot\rangle_{\bdry K}$ can be thought of as a boundary integral (for smooth enough inputs it is actually a boundary integral).
Therefore, the resulting broken ultraweak variational formulation seeks 
\begin{equation}
	\begin{gathered}
		(\fku_0,\hat{\fku})=\fku\in\scU=\scU_0\times\hat{\scU}\,,\\
		(u,\bsq)=\fku_0\in\scU_0=\Leb^2(\Omega)\times\bLeb^2(\Omega)\,,\qquad
		(\hat{u},\hat{q}_{\nml})=\hat{\fku}\in\hat{\scU}=\HSo_0^{\onehalf}(\bdry\mesh)\times\HSo^{-\onehalf}(\bdry\mesh)\,,
	\end{gathered}
	\label{eq:BrokenUWSpaces}
\end{equation}
such that
\begin{equation}
	\begin{aligned}
		&b(\fku,\fkv)=\ell(\fkv)\,\qquad\forall(v,\btau)=\fkv\in\scV=\HSo^1(\mesh)\times\bHSo(\div,\mesh)\,,\\
		&b\big((\fku_0,\hat{\fku}),\fkv\big)=b_0(\fku_0,\fkv)+\hat{b}(\hat{\fku},\fkv)\,,
			\qquad \ell\big((v,\btau)\big)=(r,v)_\mesh\,,\\
		&b_0\big((u,\bsq),(v,\btau)\big)=-(\bsq,\nabla v)_\mesh+(\textstyle{\frac{1}{k}}\bsq,\btau)_\mesh-(u,\div\btau)_\mesh\,,\\
		&\hat{b}\big((\hat{u},\hat{q}_{\nml}),(v,\btau)\big)=\langle\hat{q}_{\nml},v_{\bdry\mesh}\rangle_{\bdry\mesh}
			+\langle\hat{u},\btau_{\bdry\mesh}\rangle_{\bdry\mesh}\,,
	\end{aligned}
	\label{eq:BrokenUW}
\end{equation}
where $v_{\bdry\mesh}=\prod_{K\in\mesh}(v|_K)\big|_{\bdry K}$ and $\btau_{\bdry\mesh}=\prod_{K\in\mesh}(\btau|_K)\big|_{\bdry K}\!\cdot\!\nml_K$.
This formulation can also be proved to be well-posed, with stability properties independent of the choice of the mesh \cite{BrokenForms15,Keith16}.
With nontrivial boundary conditions, $g\neq0$, simply consider $\ell(\fkv)=(r,v)_\mesh-\langle\tilde{g}_{\bdry\mesh},\btau_{\bdry\mesh}\rangle_{\bdry\mesh}$ instead, where $\tilde{g}\in\HSo^1(\Omega)$ is an extension of $g\in\HSo^{\onehalf}(\bdry\Omega)=\{f=\tilde{f}|_{\bdry\Omega}\mid \tilde{f}\in\HSo^1(\Omega)\}$, and add $\tilde{g}$ to the solution $u$ of \eqref{eq:BrokenUW} to obtain the final temperature.

Despite looking intricate, the broken ultraweak variational formulation has the advantage of removing much of the inter-element compatibility conditions, since some of its trial variables are in $\Leb^2$ and its test variables are discontinuous along the elements.
The only inter-element compatibility is due to the skeleton variables, which reside solely on the element boundaries.
In 2D, as we mentioned before, this is extremely convenient since the element boundaries are simply 1D line segments.

\subsection{Discretization and the DPG methodology} 
\label{sec:DPGcharacteritics}

In this section we present the procedure of discretizing the ultraweak formulations.
The Galerkin method is the widely used approach for conventional formulations. It employs the same test and trial spaces, leading to a square linear system of equations.
Indeed, consider the primal formulation in \eqref{eq:Primal}, with $\{\fku^\scP_j\}_{j=1}^N$ being a basis for the discrete subspaces $\scU_h^\scP=\scV_h^\scP\subseteq\scU^\scP=\scV^\scP$. 
Then, the discrete solution $u_h=\sum_{j=1}^N(\sfu_h)_j\fku^\scP_j\in\scU_h^\scP$ for $\sfu_h\in\R^N$, satisfies
\begin{equation}
	\sfB^\scP\sfu_h=\sfell^\scP\,,
\end{equation}
where $\sfB_{ij}^\scP=b^\scP(\fku^\scP_j,\fkv^\scP_i)$ and $\sfell_i^\scP=\ell^\scP(\fkv^\scP_i)$ with $\fkv^\scP_i=\fku^\scP_i$, so that $\sfB^\scP\in\R^{N\!\times\! N}$ and $\sfell^\scP\in\R^N$.
The basis functions, $\fku^\scP_j$, are chosen with a very small support not exceeding a few neighboring elements, resulting in a computationally practical method due to the sparse structure of $\sfB^\scP$.

In general, when the trial and test spaces are different, $\scU\neq\scV$, this approach is still possible but requires finding bases $\{\fku_j\}_{j=1}^N$ and $\{\fkv_i\}_{i=1}^N$ for $\scU_h\subseteq\scU$ and $\scV_h\subseteq\scV$ respectively.
However, two issues immediately arise.
First, the canonical polynomial-based discrete basis of $\scV_h\subseteq\scV$ typically is not of size $N$ (the same size of the basis for $\scU_h$).
Second, even if a nonstandard basis for $\scV_h$ of the right size is found, the resulting numerical method could very well be unstable, meaning that the inf-sup inequality,
\begin{equation}
	\inf_{\delta\fku_h\in\scU_h\setminus\{0\}}\sup_{\fkv_h\in\scV_h\setminus\{0\}}
		\frac{b(\delta\fku_h,\fkv_h)}{\|\delta\fku_h\|_{\scU}\|\fkv_h\|_{\scV}}=\gamma_h>0\,,
	\label{eq:infsupstable}
\end{equation}
might \textit{not} hold.
In fact, depending on the equation and mesh size, even the Galerkin method can be unstable.
Minimum residual finite element methods overcome these two difficulties by design.

Let $\scU'$ and $\scV'$ be the dual spaces to $\scU$ and $\scV$ respectively, and define $\scB:\scU\to\scV'$ and its adjoint $\scB':\scV\to\scU'$ through duality pairings as $\langle\scB\fku,\fkv\rangle=b(\fku,\fkv)=\langle\fku,\scB'\fkv\rangle$.
Then, for a discrete trial space $\scU_h\subseteq\scU$, minimum residual methods seek the minimizer of the residual \cite{DPGOverview,Keith16},
\begin{equation}
	\fku_h^\opt=\argmin_{\delta\fku_h\in\scU_h}\|\scB\delta\fku_h-\ell\|_{\scV'},\quad\Leftrightarrow\quad
		b(\fku_h^\opt,\fkv^\opt)=\ell(\fkv^\opt)\quad\forall\fkv^\opt\in\scV^\opt=\scR_{\scV}^{-1}\scB\scU_h\,,
\end{equation}
where $\scR_{\scV}:\scV\to\scV'$ is the Riesz map, which is defined by duality as $\langle\scR_{\scV}\fkv,\delta\fkv\rangle=(\fkv,\delta\fkv)_\scV$, with $(\cdot,\cdot)_\scV$ being the inner product of the Hilbert space $\scV$.
Here, $\scV^\opt=\scR_{\scV}^{-1}\scB\scU_h$ is called the optimal test space, because this exact choice of discrete test space automatically results in the best inf-sup stable discrete method satisfying \eqref{eq:infsupstable} \cite{DPGOverview}.
Given an element of the basis for $\scU_h$, $\fku_i\in\{\fku_j\}_{j=1}^N$, the corresponding optimal test function is $\fkv_i^\opt=\scR_{\scV}^{-1}\scB\fku_i$. With these choices the resulting matrix $\sfB^\opt_{ij}=b(\fku_j,\fkv_i^\opt)$, called the optimal stiffness matrix, is always symmetric positive-definite.

Unfortunately, computing $\scR_{\scV}^{-1}$ is impossible since $\scV$ is infinite-dimensional.
Thus, minimum residual methods simply make a choice of an \textit{enriched} test space $\scV_r\subseteq\scV$ (with $M=\dim(\scV_r)\geq\dim(\scU_h)=N$) over which the operator is inverted.
The advantage is that this enriched space may be discretized with a standard canonical polynomial-based basis, $\{\fkv_i\}_{i=1}^M$, and ultimately the resulting \textit{near}-optimal space is $\scV_h=\scV^\nopt=\scR_{\scV_r}^{-1}\scB\scU_h$ and its corresponding \textit{near}-optimal basis is $\fkv_i^\nopt=\scR_{\scV_r}^{-1}\scB\fku_i$ for every $\fku_i\in\{\fku_j\}_{j=1}^N$.
The resulting discrete method can be shown to be equivalent to the linear system,
\begin{equation}
	\sfB^\nopt\sfu_h=\sfB^\T\sfG^{-1}\sfB\sfu_h=\sfB^\T\sfG^{-1}\sfell=\sfell^\nopt\,,
	\label{eq:DPGlinearsystem}
\end{equation}
where $\fku_h=\sum_{j=1}^N(\sfu_h)_j\fku_j\in\scU_h$ is the discrete solution; the Gram matrix $\sfG_{ij}=(\fkv_i,\fkv_j)_\scV$ is a discretization of $\scR_{\scV_r}$; $\sfB_{ij}=b(\fku_j,\fkv_i)$ and $\sfell_i=\ell(\fkv_i)$ are called the enriched stiffness matrix and load; and $\sfB^\nopt_{ij}=b(\fku_j,\fkv_i^\nopt)$ and $\sfell_i^\nopt=\ell(\fkv_i^\nopt)$ are the near-optimal stiffness matrix and load.
Clearly the enriched stiffness matrix is rectangular and tall, $\sfB\in\R^{M\!\times\! N}$ with $M\geq N$, while the near-optimal stiffness matrix is square and symmetric positive-definite, $\sfB^\nopt\in\R^{N\!\times\! N}$.
To implement, one has to form the Gram matrix ($\sfG\in\R^{M\!\times\! M}$), enriched stiffness matrix ($\sfB\in\R^{M\!\times\! N}$) and enriched load vector ($\sfell\in\R^{M}$) first, then calculate the near-optimal stiffness matrix ($\sfB^\nopt=\sfB^\T\sfG^{-1}\sfB\in\R^{N\!\times\! N}$) and near-optimal load vector ($\sfell^\nopt=\sfB^\T\sfG^{-1}\sfell\in\R^{N}$), and finally solve for the basis coefficients of the discrete solution ($\sfu_h\in\R^N$).

All this derivation holds for any arbitrary linear variational formulation including the ultraweak formulations in \eqref{eq:UW} and \eqref{eq:BrokenUW}.
The method is near-optimal in that it is designed to approximate the optimal method (with $\sfB^\opt$), so in principle it is not \textit{known} to be stable, but in practice it typically is or can be made stable (if it is not stable simply enrich $\scV_r$ even more so that $M\gg N$).
In fact, the stability of the near-optimal method can rigorously be proved by constructing a Fortin operator, $\Pi_F:\scV\to\scV_r$ \cite{QiuFortin,BrokenForms15}.

However, there are major differences between applying this method to the ultraweak formulation in \eqref{eq:UW} and the broken ultraweak formulation in \eqref{eq:BrokenUW}.
Namely, for the standard ultraweak formulation the enriched (sparse) stiffness matrix, $\sfB$, and the Gram matrix, $\sfG$, are assembled globally first and then the near-optimal stiffness matrix, $\sfB^\nopt$ is computed using \eqref{eq:DPGlinearsystem}.
This is very expensive, especially due to the inversion of $\sfG$.
Thus, despite many advantages, the method is not very practical.
On the other hand, when using broken test spaces, as in the broken ultraweak formulation, the matrix $\sfG$ has a disjoint diagonal block structure, where each block corresponds to one element.
Hence, the Gram matrix can be inverted locally, allowing the local near-optimal stiffness matrices $(\sfB^\nopt)_K$ to be computed directly for each element $K\in\mesh$.
This is turn allows $\sfB^\nopt$ to be assembled as in any other FEM.
Thus, using formulations with broken test spaces localizes the computations and parallelizes the assembly, making it a practical FEM.
However, when compared to traditional FEMs, the local computations are more expensive due to the additional skeleton variables.
Note that the broken ultraweak formulation in \eqref{eq:BrokenUW} has an enriched stiffness matrix with the structure,
\begin{equation}
	\begin{gathered}
		{}\\[-4pt]
			\text{\footnotesize$\{\fkv_i\}_{i=1}^{M}$}\left\{ \vphantom{\begin{pmatrix} | \\ | \\ | \end{pmatrix}} \right. 
			\begin{bmatrix}
				\bovermat{\{(\fku_0)_j\}_{j=1}^{N_0}}{\qquad\quad|\qquad\quad} 
					& \bovermat{\{\hat{\fku}_j\}_{j=1}^{\hat{N}}}{\quad\qquad|\qquad\quad} \\ 
				\quad\sfB_0 & \quad\hat{\sfB} \\ 
				\qquad\quad|\qquad\quad & \qquad\quad|\qquad\quad
			\end{bmatrix}
			=\sfB=
			\begin{bmatrix}
				\sfB_{uv} & \sfB_{\bsq v} & \sfB_{\hat{u}v} & \sfB_{\hat{q}_{\nml} v} \\
				\sfB_{u\btau} & \sfB_{\bsq\btau} & \sfB_{\hat{u}\btau} & \sfB_{\hat{q}_{\nml}\btau} 
			\end{bmatrix}
	\end{gathered}
\end{equation}
where $(\sfB_0)_{ij}=b_0((\fku_0)_j,\fkv_i)$ and $\hat{\sfB}_{ij}=\hat{b}(\hat{\fku}_j,\fkv_i)$, with the $\scU_h$-basis $\{\fku_j\}_{j=1}^N=\{((\fku_0)_j,0)\}_{j=1}^{N_0}\cup\{(0,\hat{\fku}_j)\}_{j=1}^{\hat{N}}$ so that $N=N_0+\hat{N}$, and similarly with the other sub-blocks.

In the literature, the application of minimum residual methods to variational formulations with broken test spaces is referred to as the \textit{DPG methodology}.
The methodology is quite general as it can be applied to variational formulations other than the broken ultraweak such as broken primal or broken mixed formulations \cite{Keith16,BrokenForms15}. 
Each application case results in a different DPG method similar to how the Galerkin methodology can be applied to primal and mixed formulations (where $\scU_h=\scV_h$).
Nonetheless, the lack of inter-element compatibility restrictions on the $\scU_0$-part of the trial space (which lies in copies of $\Leb^2$) makes the ultraweak formulation a natural candidate to develop a DPG method for polygonal elements.

It is worth mentioning that the DPG methodology carries a natural arbitrary-order residual-based \textit{a posteriori} error estimator.
The expression for the residual is,
\begin{equation}	
	\|\scB\fku_h-\ell\|_{\scV'}^2\approx\|\scR_{\scV_r}^{-1}(\scB\fku_h-\ell)\|_{\scV}^2
		=(\sfB\sfu_h-\sfell)^\T\sfG^{-1}(\sfB\sfu_h-\sfell)\,,
	\label{eq:residualexpression}
\end{equation}
where $\fku_h$ (and $\sfu_h$) is the solution.
Note that the test spaces are broken, so the computations can be performed locally. 
Therefore, \eqref{eq:residualexpression} can serve as an \textit{a posteriori} error estimator for driving different adaptive strategies \cite{Keith16,goalorientedDPG}.
Adaptivity in its own right is a very interesting subject of study for polygonal elements, as they provide great flexibility for the implementation of such strategies without resorting to constrained approximations to deal with hanging nodes. 
More details on this will be given in Section \ref{sec:results_adapt}.

A final comment on minimum residual FEMs, including all DPG methods, is that the choice of test norm (or inner product) for $\scV$, which appears in the computation of $\sfG$, has a significant influence.
Generally speaking, the standard norms are usually chosen as test norms.
For example, the standard norm for the broken ultraweak formulation in \eqref{eq:BrokenUW} is,
\begin{equation}
	\|(v,\btau)\|_{\scV}^2=\|v\|_{\HSo^1(\mesh)}^2+\|\btau\|_{\HSo(\div,\mesh)}^2
		=(v,v)_\mesh+(\nabla v,\nabla v)_\mesh+(\btau,\btau)_\mesh+(\div\btau,\div\btau)_\mesh\,.
\end{equation}
However, there are other norms that still make $\scV$ a Hilbert space but lead to different results.
Specifically for the broken ultraweak formulations, the adjoint graph norm has interesting properties \cite{DPGOverview}. 
Using the ultraweak formulation in \eqref{eq:UW}, the first two terms in this norm can be derived as,
\begin{equation}
	\|(v,\btau)\|_{\scV}^2=\|\textstyle{\frac{1}{k}}\btau-\nabla v\|_{\bLeb^2(\mesh)}^2+\|-\div\btau\|_{\Leb^2(\mesh)}^2
		+\varepsilon^2\big(\|v\|_{\Leb^2(\mesh)}^2+\|\btau\|_{\bLeb^2(\mesh)}^2\big)\,,
	\label{eq:adjointgraphnorm}
\end{equation}
where $\|\cdot\|_{\Leb^2(\mesh)}^2=(\cdot,\cdot)_\mesh$ and the same with $\|\cdot\|_{\bLeb^2(\mesh)}^2$.
The third term, which has the $\varepsilon^2$ factor, makes the norm localizable, because otherwise \eqref{eq:adjointgraphnorm} would not be a norm for arbitrary broken functions $v\in\HSo^1(\mesh)$ (although it would be a norm for $v\in\HSo_0^1(\Omega)$).
One can choose an arbitrary value for $\varepsilon>0$, but using small values of $\varepsilon$ (with the caveat of ill-conditioned local problems) is of particular interest for certain equations, such as Helmholtz \cite{GopMugaHelmholtz}.
Note that the corresponding inner products for the (real-valued) Hilbert space $\scV$ can be derived from the polarization identity, $(\fkv_1,\fkv_2)_\scV=\frac{1}{4}\big(\|\fkv_1+\fkv_2\|_\scV^2-\|\fkv_1-\fkv_2\|_\scV^2\big)$.

\subsection{Choice of trial and test spaces}
\label{sec:trialandtestspaces}

The choice of trial and test spaces is important to establish the method's convergence.
As mentioned before, strict inter-element compatibility requirements leaves very limited options.
Particularly, the problem seems to be extremely complicated for general polygons with high-order discretizations.
Fortunately, the $\scU_0$ trial space component of the broken ultraweak formulation in \eqref{eq:BrokenUWSpaces} consists of copies of $\Leb^2$, so its discretization can be discontinuous across the elements.
Moreover, the test spaces are broken, so their discretization should be discontinuous across elements too.
This freedom allows one to create bases locally, disregarding the neighboring elements.
In particular, bases may be defined by restriction (to the polygonal element of interest), as we will see next.

\begin{figure}[!ht]
 	\centering
  \includegraphics[trim=0cm 0cm 0cm 0cm,clip=true,scale=0.55]{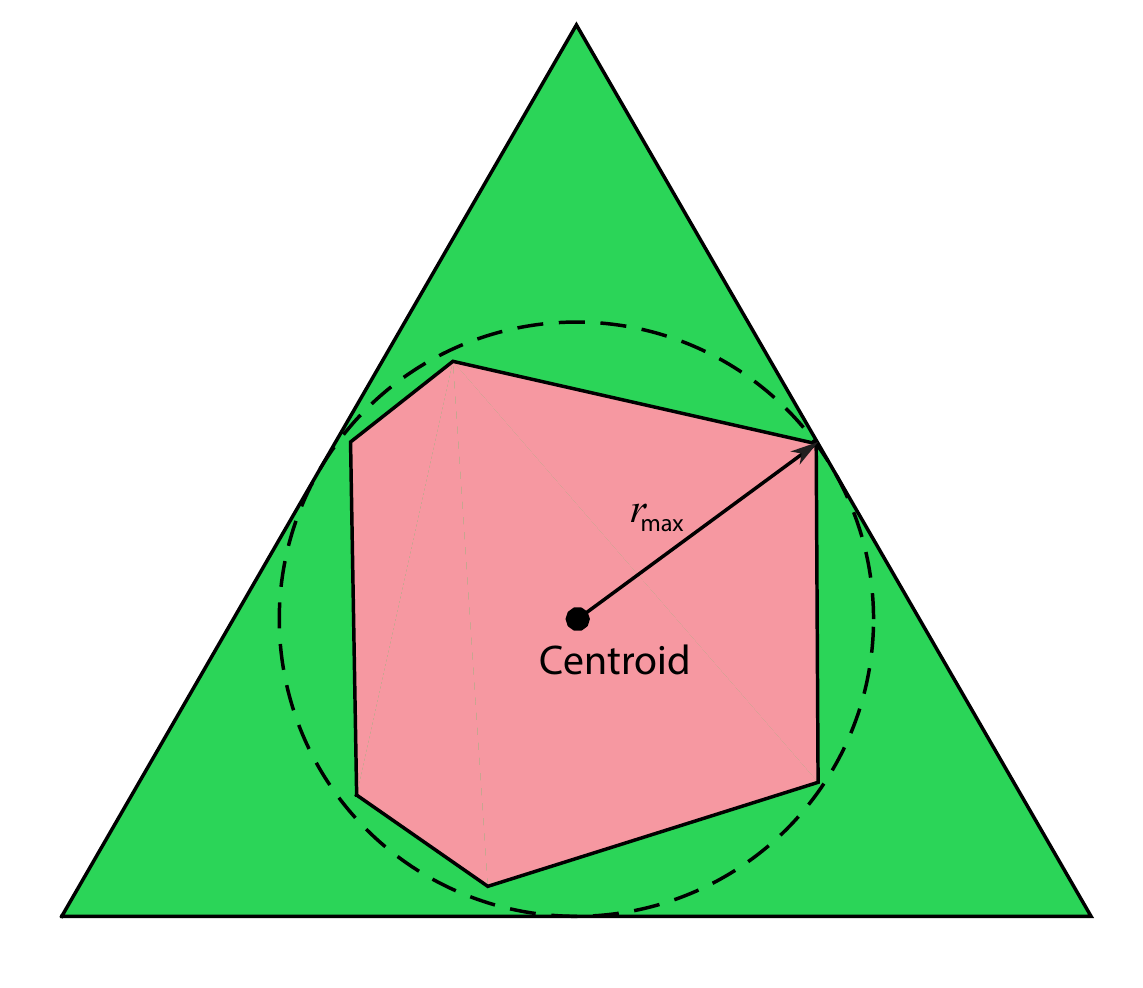}
  \caption{
	Bounding triangle of a polygonal element. The equilateral triangle is defined such that the bounding circle centered at the polygon's centroid is inscribed.
	}
  \label{fig:boundingtriangle}
\end{figure}

Our procedure is similar to that in \cite{hpdgfem} where a bounding box was utilized, but we use a bounding triangle instead.
First, the centroid of the polygon and the furthest vertex from the centroid are determined. Next, a bounding circle centered at the centroid and passing through the furthest vertex is defined.
Then, the bounding equilateral triangle inscribing the circle is computed such that one of its edge-midpoints is the polygon's furthest vertex.
This is shown in Figure \ref{fig:boundingtriangle}.
Lastly, the ``usual'' high-order polynomial shape functions for the triangle are used and then restricted to the polygon.
We use the term ``usual'' liberally, but to clarify, we include further details below.  

There are several spaces at the infinite-dimensional level which we want to discretize using this technique.
Namely, the test space components, $\HSo^1(\mesh)$ and $\bHSo(\div,\mesh)$, and the $\scU_0$ trial space component, which may be represented by $\Leb^2(\Omega)$.
Following our technique, the procedure reduces to finding the local discretizations of $\HSo^1(T_K)$, $\bHSo(\div,T_K)$ and $\Leb^2(T_K)$, where $T_K$ is the bounding triangle of the polygonal element $K\in\mesh$.
These three spaces actually form a differential de Rahm exact sequence, and it is convenient that their respective discretizations do too.
For triangles, this is satisfied by the classical N\'{e}d\'{e}lec sequence of the first type \cite{hpbook2,Fuentes2015},
\begin{equation}
	\begin{gathered}
		\xymatrix@R-1.5pc{
		\HSo^1(T_K)\ar[r]^{\curl\quad} &\bHSo(\div,T_K) \ar[r]^{\nabla\cdot}&\,\,\Leb^2(T_K)\,\,\\
		\mcP^p(T_K) \ar@{}[u]|*=0[@]{\subseteq} \ar[r]^{\curl\quad} 
			&\mcR\mcT^p(T_K) \ar@{}[u]|*=0[@]{\subseteq} \ar[r]^{\nabla\cdot} 
				&\mcP^{p-1}(T_K) \ar@{}[u]|*=0[@]{\subseteq}\,,
		}
	\end{gathered}
	\label{eq:exactsequence}
\end{equation}
where $\mcP^p(T_K)$ are the polynomials in $x=(x_1,x_2)$ of total order less than or equal $p\in\N$, the 2D Raviart-Thomas space is $\mcR\mcT^p(T_K)=(\mathcal{P}^{p-1}(T_K))^2+x\mathcal{P}^{p-1}(T_K)$  (a rotation of the 2D N\'{e}d\'{e}lec space), and the 2D scalar-to-vector curl operator is defined as $\curl(u)=\big(\begin{smallmatrix}0&1\\-1&0\end{smallmatrix}\big)\nabla u$ for any $u\in\HSo^1(T_K)$.
Notice that the parameter $p$ represents the order of the discrete \textit{sequence} and does not necessarily coincide with the order of the polynomials of a particular discretization.
For example if $p=3$, the discretization of $\Leb^2(T_K)$ are the polynomials of at most total order $p-1=2$.

This sequence has many desirable properties, and precisely because of these, we prefer to use a bounding triangle instead of a bounding box.
In particular, the spaces are invariant under affine transformations (the spaces remain the same even if the bounding triangle is arbitrarily rotated about the polygon centroid); the overall drop of polynomial order across the sequence is one (from $\mcP^p(T_K)$ to $\mcP^{p-1}(T_K)$); the approximation properties are suitable (see Appendix \ref{sec:AppConvergence}); and they are the smallest possible spaces with all these properties (see \cite[\S3.4]{arnold2006finite}).

Having said that, a similar procedure can be carried out for a bounding box, $Q_K$ of $K\in\mesh$, where the spaces become 
\begin{equation}
	\begin{gathered}
		\xymatrix@R-1.5pc{
		\HSo^1(Q_K)\ar[r]^{\curl\qquad\qquad} &\qquad\bHSo(\div,Q_K)\qquad \ar[r]^{\qquad\nabla\cdot}&\,\,\Leb^2(Q_K)\,\,\\
		\mcQ^{p,p}(Q_K) \ar@{}[u]|*=0[@]{\subseteq} \ar[r]^{\curl\qquad\qquad} 
			&\mcQ^{p,p-1}(Q_K)\!\times\!\mcQ^{p-1,p}(Q_K) \ar@{}[u]|*=0[@]{\subseteq} \ar[r]^{\qquad\nabla\cdot} 
				&\mcQ^{p-1,p-1}(Q_K) \ar@{}[u]|*=0[@]{\subseteq}\,,
		}
	\end{gathered}
	\label{eq:exactsequencequad}
\end{equation}
with $\mcQ^{p,q}(Q_K)=\mcP^p(x_1)\otimes\mcP^q(x_2)$.

In either case, the final spaces for the polygon $K\subseteq T_K$ (or $K\subseteq Q_K$) are defined by restricting the domain to $K\in\mesh$, so we denote them by $\mcP^p(K)$ and $\mcR\mcT^p(K)$ (or $\mcQ^{p,p}(K)$) instead.

The only remaining spaces to specify are those of the skeleton variables lying in the $\hat{\scU}$ trial space component (see \eqref{eq:BrokenUWSpaces}).
These can also be deduced using the same philosophy of exact sequences, but utilizing the traces instead.
Indeed, the spaces $\HSo_0^{\onehalf}(\bdry\mesh)$ and $\HSo^{-\onehalf}(\bdry\mesh)$ are merely $\mesh$-tuples of compatible traces of $\HSo^1(K)$ and normal-traces of $\bHSo(\div,K)$ respectively.
If two elements of different type (a triangle and a quadrilateral) share an edge, the discrete spaces should be compatible across that edge.
This is the case when considering the $\HSo^1$-discretizations of triangles and quadrilaterals: even though the discretizations themselves are different ($\mcP^p$ and $\mcQ^{p,p}$), their restrictions to edges are exactly the same, $\mcP^p(e)$, where $e$ represents an edge parametrized linearly by $t_e$.
The same occurs with the $\bHSo(\div)$-discretizations, which have $\mcP^{p-1}(e)$ as normal-trace along the edges. 
Additionally, the $\HSo^1$-discretizations should be compatible at vertices.
This is consistent with 1D discretizations of $\HSo^1$ and $\Leb^2$, which also form an exact sequence, but instead occurring along the boundary of each element and being edge-parametrized along all edges (see \cite[\S1.6]{Fuentes2015}).
This pattern should hold for arbitrary polygons as well.
For this, let $\mcE(K)$ be the set of edges of a polygon $K\in\mesh$, and define the local discretizations,
\begin{equation}
	\begin{aligned}
		\mcP^{p-1}(\bdry K)&=\{\hat{w}_K\mid\hat{w}_K|_{e}\in\mcP^{p-1}(e),\,\forall e\in\mcE(K)\}
			\subseteq\HSo^{-\onehalf}(\bdry K)\,,\\
		\mcP^{p}_C(\bdry K)&=\mcP^p(\bdry K)\cap C^0(\bdry K)\subseteq\HSo^{\onehalf}(\bdry K)\,,
	\end{aligned}
	\label{eq:localtracediscretizations}
\end{equation}
where $C^0(\bdry K)$ are the continuous functions in $\bdry K$ (the intersection ensures that values of neighboring edges coincide at a common vertex), and the local trace spaces are $\HSo^{\onehalf}(\bdry K)=\{\hat{u}_K=u|_{\bdry K}\mid u\in\HSo^1(K)\}$ and $\HSo^{-\onehalf}(\bdry K)=\{(\hat{q}_{\nml})_K=\bsq|_{\bdry K}\!\cdot\!\nml_K\mid \bsq\in\bHSo(\div,K)\}$.

Now we have enough information to actually globally define the discrete trial space.
For a value of $p\in\N$, it is
\begin{equation}
	\begin{aligned}
		\scU_h=\big\{(u,\bsq,\hat{u},\hat{q}_{\nml})\in\scU\mid u|_K&\in\mcP^{p-1}(K),\,\bsq|_K\in\big(\mcP^{p-1}(K)\big)^2,\\
			&\qquad\hat{u}_K\in\mcP^{p}_C(\bdry K),\,(\hat{q}_{\nml})_K\in\mcP^{p-1}(\bdry K),\,\forall K\in\mesh\big\}\,.
	\end{aligned}
	\label{eq:Uhdefinition}
\end{equation} 
Notice that the condition $(u,\bsq,\hat{u},\hat{q}_{\nml})\in\scU$ (so $(\hat{u},\hat{q}_{\nml})\in\hat{\scU}$) implies that $\hat{u}$ vanishes at the boundaries, that $\hat{u}_{K_1}|_e=\hat{u}_{K_2}|_e$, and that $(\hat{q}_{\nml})_{K_1}|_e=-(\hat{q}_{\nml})_{K_2}|_e$, where $e$ is a common edge between the elements $K_1$ and $K_2$.
No such compatibility implications exist for $(u,\bsq)\in\scU_0$.

For the enriched test space, the discretizations are chosen from a sequence of order $p+\Delta p$, and we say the space is $p$-enriched, so that
\begin{equation}
	\scV_r=\big\{(v,\btau)\mid v|_K\in\mcP^{p+\Delta p_K}(K),\,\btau|_K\in\mcR\mcT^{p+\Delta p_K}(K),\,\forall K\in\mesh\big\}\,.
	\label{eq:Vrdefinition}
\end{equation}
The notation $\Delta p_K$ indicates that this value is element-dependent.
In fact, recall that for minimum residual methods to work, $M=\dim(\scV_r)\geq\dim(\scU_h)=N$, and this restriction on the dimensionality should hold locally as well.
Thus, $\Delta p_K$ has to be chosen such that this condition holds.
This is important for the polygonal element methods, because when a polygon has many sides, the size of the local trial space may be quite large and a large value of $\Delta p_K$ must be chosen for that particular element.

To elaborate, consider an interior $n$-sided polygonal element $K$ (so that $\bdry K\cap\bdry\Omega=\varnothing$).
Its local trial and test space dimensions would be
\begin{equation}
	\begin{aligned}
		\dim(\scU_h(K))&=\overbrace{\textstyle{\frac{1}{2}}p(p+1)}^{u|_K}+\overbrace{p(p+1)}^{\bsq|_K}
			+\overbrace{n+n(p-1)}^{\hat{u}_K}+\overbrace{np}^{(\hat{q}_{\nml})_K}\,,\\
		\dim(\scV_r(K))&=\underbrace{\textstyle{\frac{1}{2}}(p+\Delta p_K+1)(p+\Delta p_K+2)}_{v|_K}
			+\underbrace{(p+\Delta p_K)(p+\Delta p_K+2)}_{\btau|_K}\,.
	\end{aligned}
	\label{eq:deltapadhoc}
\end{equation}
Thus, for $p=2$ and $n=3$ (a triangle), $\dim(\scU_h(K))=21$, so that a value of $\Delta p_K=1$ is sufficient ($\dim(\scV_r(K))=25$); but if $p=2$ and $n=8$ (an octagon), $\dim(\scU_h(K))=41$, a value of at least $\Delta p_K=3$ (so that $\dim(\scV_r(K))=56$) is required.
Having said that, sometimes for simplicity a valid value of $\Delta p$ is chosen uniformly throughout the mesh (this is the case for all of our examples in Section \ref{sec:results}). 

\begin{figure}[htp]
	\centering
	\includegraphics[width=15.5cm]{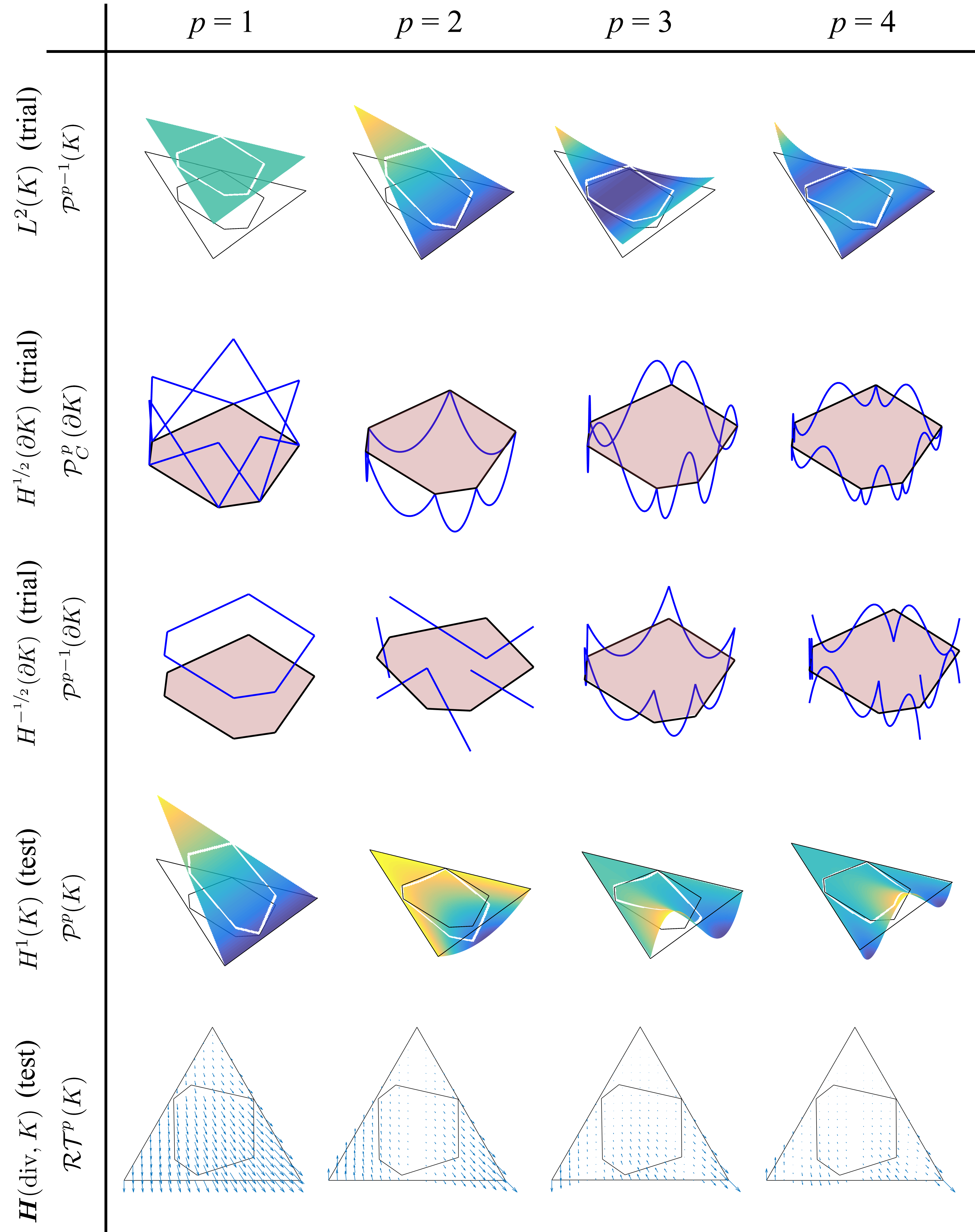}
	\caption{Some of the shape functions on a polygonal element used either as trial or test variables in the PolyDPG method. They are classified by the energy space ($\HSo^1(K)$, $\bHSo(\div,K)$, $\Leb^2(K)$ and their respective traces) and shown for different values of the parameter $p\in\N$ denoting the order of the differential \textit{sequence}. The underlying hierarchical shape functions for the bounding triangle and edges are taken from \cite{Fuentes2015}.}
	\label{fig:shape_functions}
\end{figure}

To illustrate, some representative shape functions of the components of $\scU_h(K)$ and $\scV_r(K)$ are shown in Figure \ref{fig:shape_functions} for the different energy spaces and multiple values of $p$.

We refer to the high-order polygonal DPG method resulting from this choice of trial and enriched test spaces as a PolyDPG method for Poisson's equation.
However, it can easily be generalized to ultraweak formulations coming from other linear equations (see Remark \ref{rmk:generalequations} later), so it is more appropriate to allude to a \textit{family} of PolyDPG methods. 
Note that the methods seem to be very expensive due to the large number of variables in the trial space $\scU_h$, but this is deceiving.
In fact, all of the $\scU_0$ trial space components can be statically condensed locally for ultraweak formulations, meaning that this part of the near-optimal stiffness matrix, $\sfB^\nopt$, can be effectively removed by taking Schur complements. Thus, the only remaining connectivity is that coming from the skeleton variables in $\hat{\scU}$.
So computationally speaking, solving with these variational formulations is not as costly as one might initially imagine.

\subsection{Convergence}

Since the subspaces used to discretize the ultraweak variational formulation are, rigorously speaking, subsets of the infinite dimensional trial and test spaces, PolyDPG methods are conforming FEMs.
Thus, the ``standard'' convergence theory can be applied.
However, this is an understatement because the skeleton variables are not standard, so they require a careful treatment.
The details are left to Appendix \ref{sec:AppConvergence}, but the main result is stated here along with the key assumptions.

\begin{definition*}
A collection of subsets of $\R^2$, $\mesh_{\mcK}$, is said to have the finite overlap condition if
\begin{equation}
	\overlap(\mesh_{\mcK})=\sup_{x\in\R^2}\overlap(x)<\infty\,,\qquad\overlap(x)=|\{\mcK\in\mesh_{\mcK}\mid x\in\mcK\}|\,.
\end{equation}
For a family of such collections given by a parameter $\fkh\in\fkH$, $\{\mesh_{\mcK,\fkh}\}_{\fkh\in\fkH}$, the finite overlap condition is said to be robust in $\fkh$ if there exists an integer $M_{\overlap}>0$, independent of $\fkh$, such that $\overlap(\mesh_{\mcK,\fkh})\leq M_{\overlap}$ for any $\fkh\in\fkH$.
\end{definition*}

\begin{definition*}
A triangulation $\scT(K)=\{\scT_i(K)\}_{i\in I_K}$ (with $I_K$ finite) of a (simple) polygonal element $K$ is said to be edge-compatible if for each edge of $K$, only one $\scT_i(K)$ shares that edge.
For any polygon such a triangulation is known to exist \cite{polygonears,polytriangulation,randompolytriangulation}.
The triangulation is additionally said to be shape-regular if all $\scT_i(K)$ satisfy a kind of uniform shape-regularity condition (e.g. they satisfy a minimum angle condition or the ratio of their diameters to their incircle radii remains bounded). 
\end{definition*}

\begin{theorem}
\label{thm:convergence}
Let $p\in\N$ be a polynomial order and $\{\mesh_\fkh\}_{\fkh\in\fkH}$ be a family of polygonal meshes discretizing the domain $\Omega$, such that there exist shape-regular edge-compatible triangulations for all $K\in\mesh_\fkh$ with a robust shape-regularity condition independent of $K\in\mesh_\fkh$ across all $\fkh\in\fkH$.
Assume that the associated collections of bounding triangles \textup{(}see Section \ref{sec:trialandtestspaces}\textup{)}, $\{\mesh_{T,\fkh}\}_{\fkh\in\fkH}=\{\{T_K\}_{K\in\mesh_\fkh}\}_{\fkh\in\fkH}$, where $T_K$ is the bounding triangle of a polygonal element $K$, satisfy a robust finite overlap condition.
Also, assume the existence of a linear and continuous Fortin operator, $\Pi_F:\scV\to\scV_r$, satisfying the orthogonality condition, $b(\fku_h,\fkv-\Pi_F\fkv)=0$, for all $\fku_h\in\scU_h$ and $\fkv\in\scV$, and with a continuity bound, $M_F>0$ \textup{(}so $\|\Pi_F\fkv\|_{\scV}\leq M_F\|\fkv\|_{\scV}$ for all $\fkv\in\scV$\textup{)}, where $b:\scU\times\scV\to\R$, $\ell:\scV\to\R$, $\scU_h$ and $\scV_r$ are given in \eqref{eq:BrokenUWSpaces}, \eqref{eq:BrokenUW}, \eqref{eq:Uhdefinition} and \eqref{eq:Vrdefinition}.
Then, the problem of finding $\fku_h\in\scU_h$ such that
\begin{equation*}
	b(\fku_h,\fkv_h)=\ell(\fkv_h)\,,\qquad\forall\fkv_h\in\scV_h=\scR_{\scV_r}^{-1}\scB\scU_h\,,
\end{equation*}
has a unique solution.
When compared to the unique solution of the infinite dimensional problem, $\fku\in\scU$ \textup{(}so $b(\fku,\fkv)=\ell(\fkv)$ for all $\fkv\in\scV$\textup{)}, and assuming it is regular enough, $\fku\in\scU^s\subseteq\scU$, for an $s>\frac{1}{2}$, the following $h$-convergence estimate holds provided $M_F$ is independent of $\fkh$,
\begin{equation*}
	\|\fku-\fku_h\|_{\scU}\leq Ch^{\min\{s,p\}}\|\fku\|_{\scU^s}\,,
\end{equation*}
where $h=\sup_{K\in\mesh_\fkh}\diam(K)$ and $C=C(s,p,\Omega)>0$ is a constant independent of $\fkh$ \textup{(}and so of $h$ as well\textup{)}.
For more details about $s>\frac{1}{2}$ and $\scU^s$ see Appendix \ref{sec:AppConvergence}.
Moreover, if $M_F$ is $p$-independent as well, then in the $p$-asymptotic limit $C=\tilde{C}(\ln p)^2p^{-s}$ where $\tilde{C}=\tilde{C}(s,\Omega)$ is independent of $p$.
\end{theorem}

\begin{remark}
The robust finite overlap condition is also assumed in \cite{hpdgfem}, and is not a very restrictive assumption.
It is used in the proof to establish a robust finite constant for the global $\Leb^2(\Omega)$ convergence estimates (details are in Appendix \ref{sec:AppConvergence}). 
On the other hand, the robust shape-regular edge-compatible triangulation of all elements is a more restrictive assumption, but it is necessary to prove the convergence estimates of the skeleton variables. 
\end{remark}

\begin{remark}
\label{rmk:generalequations}
As shown in Appendix A, the theorem actually holds for any well-posed broken ultraweak variational formulation with trial variables in $\Leb^2(\Omega)$ and skeleton (also trial) variables in subsets of $\HSo^{\onehalf}(\bdry\mesh)$ and $\HSo^{-\onehalf}(\bdry\mesh)$.
Thus, this result also holds for other equations such as linear elasticity, acoustics, and convection-dominated diffusion.
\end{remark}

\begin{remark}
The arguments can be easily extended to a 3D mesh with polyhedral elements provided all the faces of the polyhedra are triangular.
Then, the proof would even hold for equations involving skeleton variables representing the traces of $\bHSo(\curl,\Omega)$ spaces, like an ultraweak formulation of Maxwell's equations (see \cite{BrokenForms15}).
However, the problem (and the corresponding numerical implementation) is more challenging for general polyhedra in 3D.
\end{remark}

%




%% file: 3_Results.tex
\section{Numerical examples} 
\label{sec:results}

In this section we consider several examples to examine the performance of the PolyDPG method.
In all cases, Poisson's equation representing the nondimensionalized steady-state heat equation was solved in the domain $\Omega=(0,1)^2$.
Unless otherwise stated, bounding triangles were utilized (as opposed to bounding boxes) and the (nondimensional) conductivity was taken as $k=1$.
Also, a default uniform value of $\Delta p=1$ was used, but was increased (uniformly across the mesh, for the sake of simplicity) if deemed necessary (see \eqref{eq:deltapadhoc} in Section \ref{sec:trialandtestspaces}).
For all computations, the adjoint graph norm written in \eqref{eq:adjointgraphnorm} with $\varepsilon=1$ was used as the test space norm.

In the first example, we studied nontrivial meshes with $n$-sided convex polygons. 
In the second example, we considered highly distorted concave elements in the mesh.
The third example was inspired by problems in geoscience, where arbitrary faults separating different material properties occur.
To model this, we cut a uniform grid at an angle, so that the resulting mesh had different polygons (pentagons, quadrilaterals and triangles) with discontinuous material properties at each side of the cut.
In these three examples, ``uniform'' refinements were analyzed for different values of $p\in\N$, in the sense that the largest element diameter was roughly cut in half with each refinement.
In the final example, we described a polygonal adaptivity scheme by using the PolyDPG arbitrary-order \textit{a posteriori} error estimator, and compared it with conventional adaptive methods (using standard element shapes).
This is particularly important since adaptive refinement algorithms applied to polygonal elements have applications in topology optimization \cite{pfem_optimization,vem_optimization,antonietti2017,polytopmatlab} and crack propagation \cite{pfem_cracks,fracturemeshbias}. 

Note that in all examples we only report the relative error in the $\scU_0$ trial space component. This is because a rigorous computation of the norms in the $\hat{\scU}$ trial space component is simply not viable.
The $\scU_0$ relative error is defined as
\begin{equation}
	\label{eq:results_error}
	\text{Relative error}=\frac{\|\fku_0-(\fku_0)_h\|_{\scU_0}}{\|\fku_0\|_{\scU_0}}\,,\quad\text{with}\quad
		\|(u,\bsq)\|_{\scU_0}^2=\|u\|_{\Leb^2(\Omega)}^2+\|\bsq\|_{\bLeb^2(\Omega)}^2=(u,u)_\Omega+(\bsq,\bsq)_\Omega\,,
\end{equation}
where $\fku_0$ is the exact solution and $(\fku_0)_h$ is the computed solution from the PolyDPG method.


\begin{remark}[\textbf{\texttt{PolyDPG} software}]
Implementation of PolyDPG methods may deceptively appear difficult when compared to typical FEM algorithms, so we developed an open-source code written in MATLAB\textsuperscript{\textregistered} also called \texttt{PolyDPG} \cite{polydpg}.
It can be run sequentially or in parallel, and it supports both conventional and polygonal elements.
We hope this removes some qualms related to the implementation and makes DPG methods more accessible to other researchers. 
The shape functions used in the code were originally described in \cite{Fuentes2015} (see Figure \ref{fig:shape_functions}). 
The numerical integration was carried out by splitting the polygons into triangles (through Delaunay triangulation), followed by using Gaussian quadrature for each triangle (the Gaussian quadrature points and weights were carefully mapped back from a square), so that polynomial integrands of a certain order were computed up to machine precision. 
\end{remark}

\subsection{Mesh with convex polygons} 
\label{sec:results_poly}

\begin{figure}[!ht]  
	\centering
	\includegraphics[width=14cm]{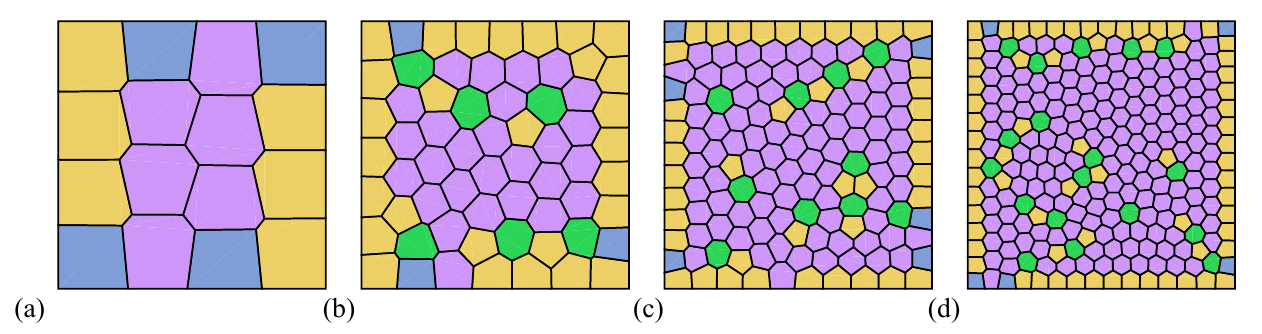}
	\caption{Four refinements of a mesh with $n$-sided convex polygonal elements. The elements are colored according to their number of sides.}
	\label{fig:poly_mesh}
\end{figure}

In this example, we investigated meshes with $n$-sided convex polygonal elements. 
The software \texttt{PolyMesher} \cite{polymesher} was used to generate the polygonal meshes.
In Figure \ref{fig:poly_mesh} an initial mesh and three subsequent refinements are displayed. 
The elements are colored according to their number of sides, ranging from $4$ (quadrilaterals) to $7$ (heptagons).
We used the manufactured solution,
\begin{equation}
	u(x,y)=\sin (\pi x)\sin (\pi y)\,,
	\label{eq:manuf1}
\end{equation}
for $(x,y)\in\Omega=(0,1)^2$ to determine the forcing, i.e. the internal heat source $r$ in \eqref{eq:poissonheat}, and the boundary conditions of $u$ at $\partial\Omega$.

\begin{figure}[!ht]
	\centering
	\includegraphics[width=16cm]{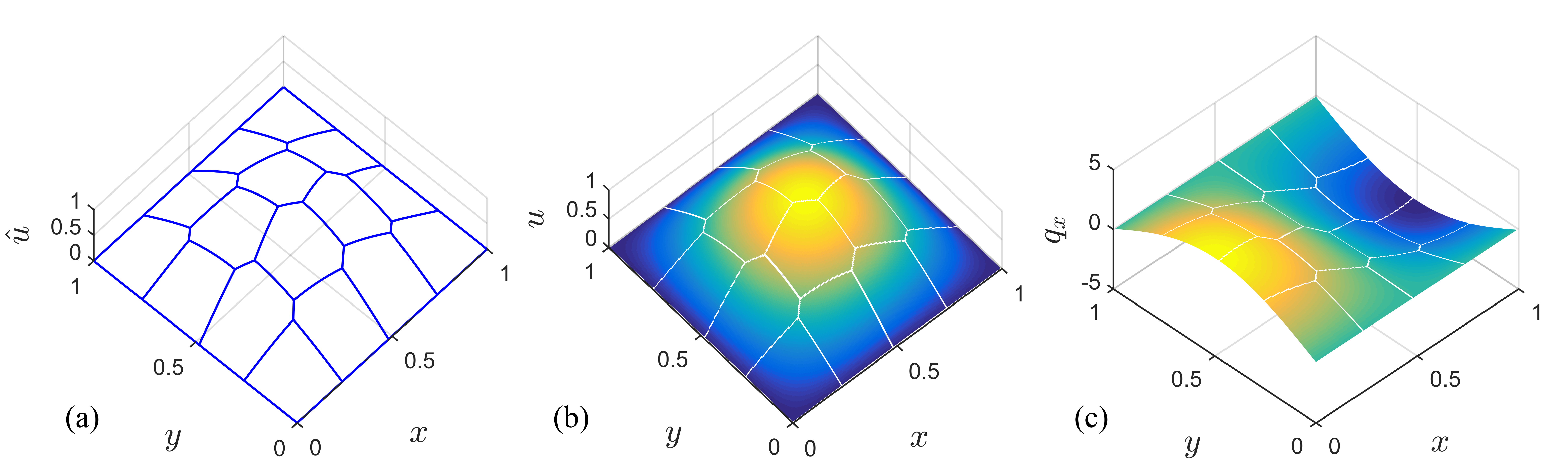}
	\caption{Numerical results using the manufactured solution in \eqref{eq:manuf1} on the coarse mesh from Figure \ref{fig:poly_mesh}(a) using $p=4$ and $\Delta p=3$: (a) skeleton temperature, (b) temperature, (c) first component of the heat flux.}
	\label{fig:poly_results}
\end{figure}

\begin{figure}[!ht]
	\centering
	\includegraphics[width=8.5cm]{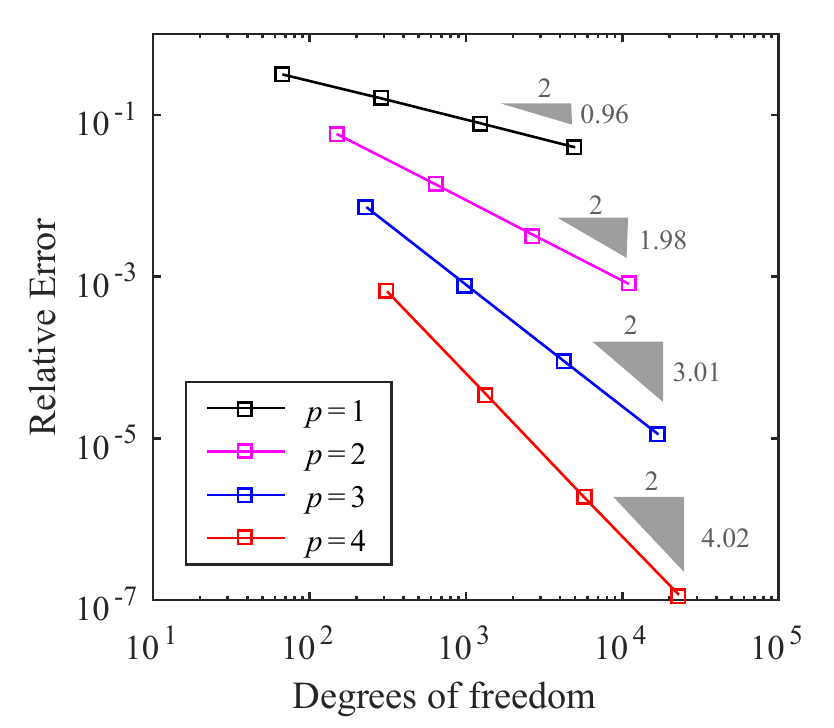}
	\caption{Convergence study of the PolyDPG method in terms of degrees of freedom. The $h$-convergence behavior is displayed for different values of $p$ using the polygonal meshes in Figure \ref{fig:poly_mesh}.}
	\label{fig:poly_error}
\end{figure}

As mentioned before, given a trial space associated to a parameter $p$, the corresponding (uniform) value of $\Delta p$ was calculated from \eqref{eq:deltapadhoc} (using the polygon with the greatest number of sides).
Given the presence of hexagons and heptagons, this meant that $\Delta p=2$ was required when $p=1,2$, while $\Delta p=3$ was needed when $p=3,4$. 
The numerical results are plotted and presented in Figure \ref{fig:poly_results} for $p=4$, including the skeleton temperature, temperature, and heat flux.
Additionally, the relative error, calculated using \eqref{eq:results_error}, is shown in Figure \ref{fig:poly_error}, where the expected $h$-convergence rates can be observed for all values of $p$ (the behavior is of the form $h^p$ as established by Theorem \ref{thm:convergence}).
Note that the number of degrees of freedom, $N_{\mathrm{dof}}$, is proportional to $h^2$. Thus, the log-log slope indicators in Figure \ref{fig:poly_error} display a $2$ in the $N_{\mathrm{dof}}$-direction, while the other label corresponds to the $h$-convergence rate, $\tilde{p}$ (so that $\frac{\tilde{p}}{2}$ is the $N_{\mathrm{dof}}$-convergence rate).

\subsection{Mesh with distorted elements} 
\label{sec:results_dist}

\begin{figure}[!ht]
	\centering
	\includegraphics[width=14cm]{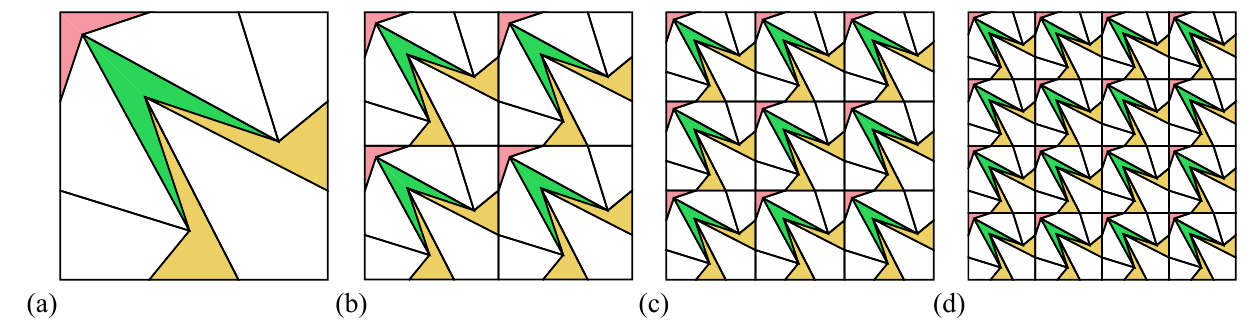}
	\caption{Four refinements using the tessellation of a mesh with highly distorted quadrilaterals. The concave elements are colored.}
	\label{fig:dist_mesh}
\end{figure}

\begin{figure}[!ht]
	\centering
	\includegraphics[width=16cm]{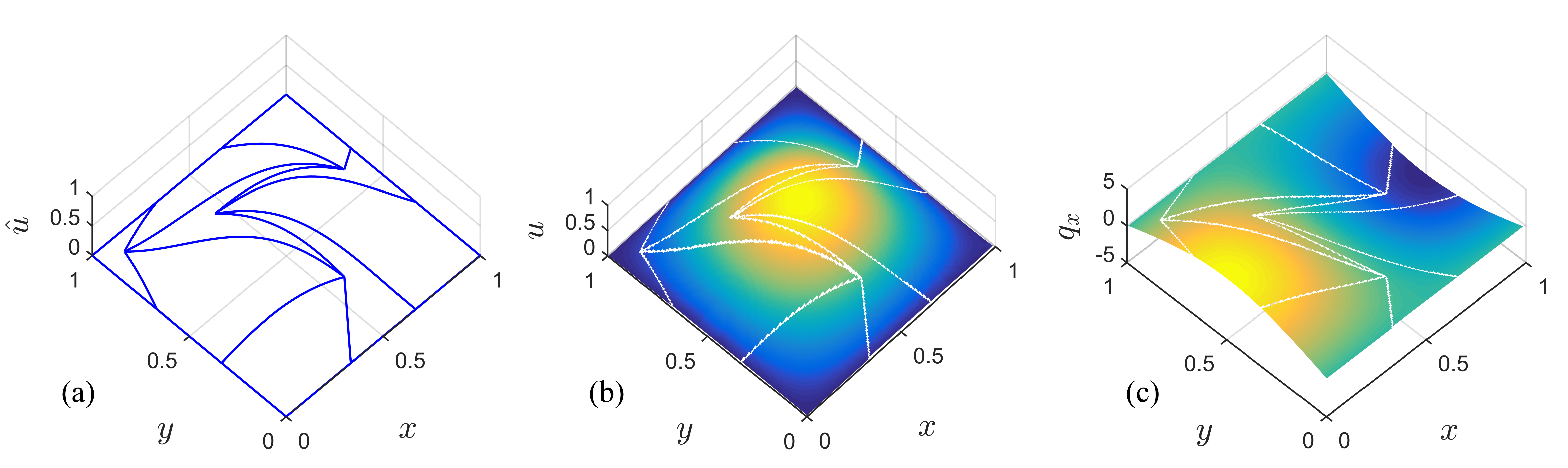}
	\caption{Numerical results using the manufactured solution in \eqref{eq:manuf1} on the coarse mesh from Figure \ref{fig:dist_mesh}(a) using $p=4$ and $\Delta p=2$: (a) skeleton temperature, (b) temperature, (c) first component of the heat flux.}
	\label{fig:dist_plots}
\end{figure}

\begin{figure}[!ht]
	\centering
	\includegraphics[width=8.5cm]{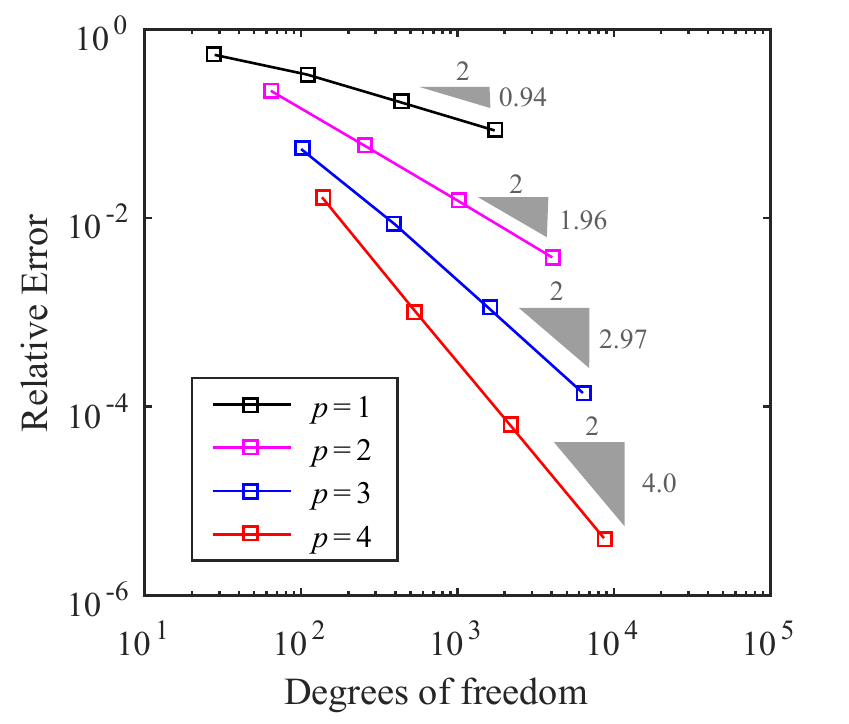}
	\caption{Convergence study of the PolyDPG method in terms of degrees of freedom. The $h$-convergence behavior is displayed for different $p$ using the meshes with highly distorted quadrilaterals in Figure \ref{fig:dist_mesh}.}
	\label{fig:dist_error}
\end{figure}

To demonstrate the distortion tolerance of PolyDPG methods, we considered a mesh with highly distorted quadrilaterals, including concave elements.
The pattern was then scaled and tessellated to produce the refinements shown in Figure \ref{fig:dist_mesh}.
This example is challenging in the sense that other numerical methods likely fail due to the degeneration of either the parametric mapping or the barycentric coordinates associated with the highly distorted elements \cite{parametricdistortion,concavebarycentric}.
The same problem as in Section \ref{sec:results_poly} was solved (see \eqref{eq:manuf1} for manufactured solution). 
The solution values and $h$-convergence rates for $1\leq p\leq4$ are shown in Figures \ref{fig:dist_plots} and \ref{fig:dist_error} respectively. 
The expected convergence behavior was observed, showing the flexibility of PolyDPG methods to deal with irregular elements.

\subsection{Interface problem} 
\label{sec:results_fault}

\begin{figure}[!ht]
	\centering
	\includegraphics[width=5cm]{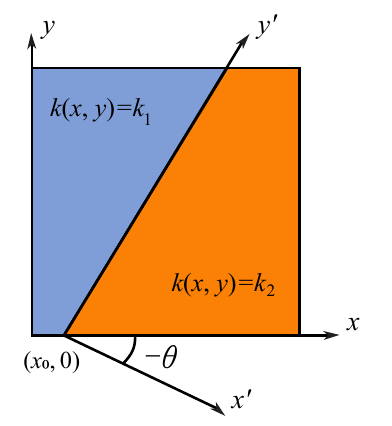}
	\caption{Material properties and rotated coordinates of the interface problem.}
	\label{fig:fault_prob}
\end{figure}

The inspiration behind this example came from geoscience applications where faults abruptly separate the material properties within a domain.
Here we considered a domain composed of two materials with different heat conductivities, which share an interface (for simplicity a straight line at an arbitrary angle dividing the square). 
The heat conductivities are assumed to be uniform on each side of the interface, taking values $k_1$ and $k_2$, as depicted in Figure \ref{fig:fault_prob}. 

\begin{figure}[!ht]
	\centering
	\includegraphics[width=14cm]{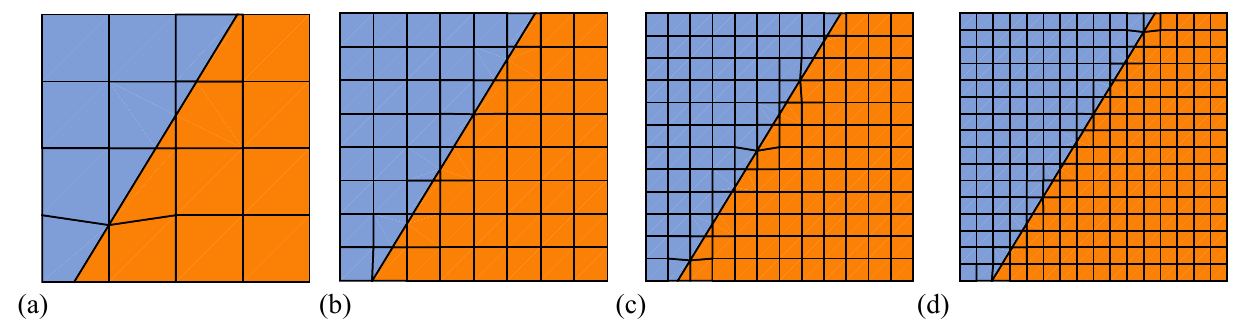}
	\caption{Four refinements of a mesh with an interface between two materials. Notice that some nodes are collapsed to a node on the interface. This is due to eliminating the undesired tiny triangles that cause ill conditioning.}
	\label{fig:fault_mesh}
\end{figure}

To model certain interfaces one would need unstructured grids. 
However, by using PolyDPG methods we are able to consider a uniform background grid and simply cut the elements through the interface, leading to the creation of triangles, right trapezoids and pentagons near the interface. 
In fact, to refine the mesh, first the background mesh was uniformly refined, and then the elements were cut by the interface line.
There is one caveat which is only evident for high values of $p$ or small values of $h$: when extremely small triangles (compared to their neighbors) are formed, the assembled stiffness matrix becomes ill conditioned (so the \textit{infinite}-precision result in Theorem \ref{thm:convergence} seizes to hold).
Thus, it is necessary to either relocate the nodes along the interface or to collapse the nodes of the small triangle into a single node on the interface.
We chose to implement the latter approach whenever the area of the small triangle was less than 1\% of the area of the background grid elements.
The meshes obtained are shown in Figure \ref{fig:fault_mesh}. 

For this problem we designed a manufactured solution that guarantees continuity of the temperature and the heat flux across the interface, taking into account the finite jump in the conductivity coefficient. 
By means of a translated and rotated system of coordinates, and following the notation in Figure \ref{fig:fault_prob}, the exact solution is given by,
\begin{equation}
	\label{eq:manuf3}
	\centering
	u(x',y') = \left\{\begin{array}{lr}
	        k_2\sin(\pi x')\sin(\pi y'), & \text{for } x'\leq 0\,,\\
	        k_1\sin(\pi x')\sin(\pi y'), & \text{for } x' > 0\,,\\
	        \end{array} \right.
\end{equation}
where the coordinates $x'$ and $y'$ come from a translation and rotation of the reference system defined by the following transformation,
\begin{equation}
	\label{eq:coord}
	\centering
	\left(\begin{array}{l}
	x' \\
	y' \\
	\end{array}\right)
	 = \left(\begin{array}{rl}
	         \cos \theta & \sin \theta\\
	        -\sin \theta & \cos \theta\\
	        \end{array} \right)
	    \left(\begin{array}{c}
	        x-x_0\\
	        y\\
	        \end{array} \right)\ .   
\end{equation}
The values of conductivity and the geometric data used for the numerical computation are $k_1=1$, $k_2=5$, $x_0=0.12$ and $\theta=\tan^{-1}(1/0.65)$. 
The nonzero boundary conditions were imposed using projection-based interpolation of the manufactured solution on the boundary edges \cite{demkowicz2008polynomial,hpbook2}.
 
\begin{figure}[!ht]
	\centering
	\includegraphics[width=16cm]{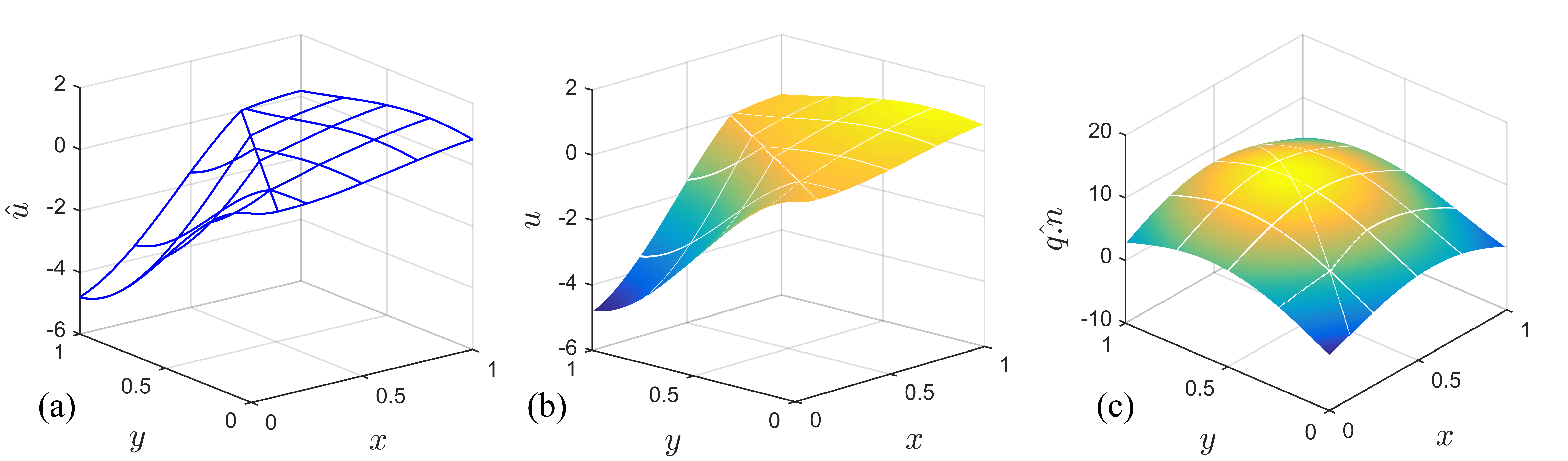}
	\caption{Numerical results using the manufactured solution in \eqref{eq:manuf3} and \eqref{eq:coord} on the coarse mesh from Figure \ref{fig:fault_mesh}(a) using $p=4$ and $\Delta p=2$: (a) skeleton temperature, (b) temperature, (c) first component of the heat flux.}
	\label{fig:fault_plots}
\end{figure}

\begin{figure}[!ht]
	\centering
	\includegraphics[width=8.5cm]{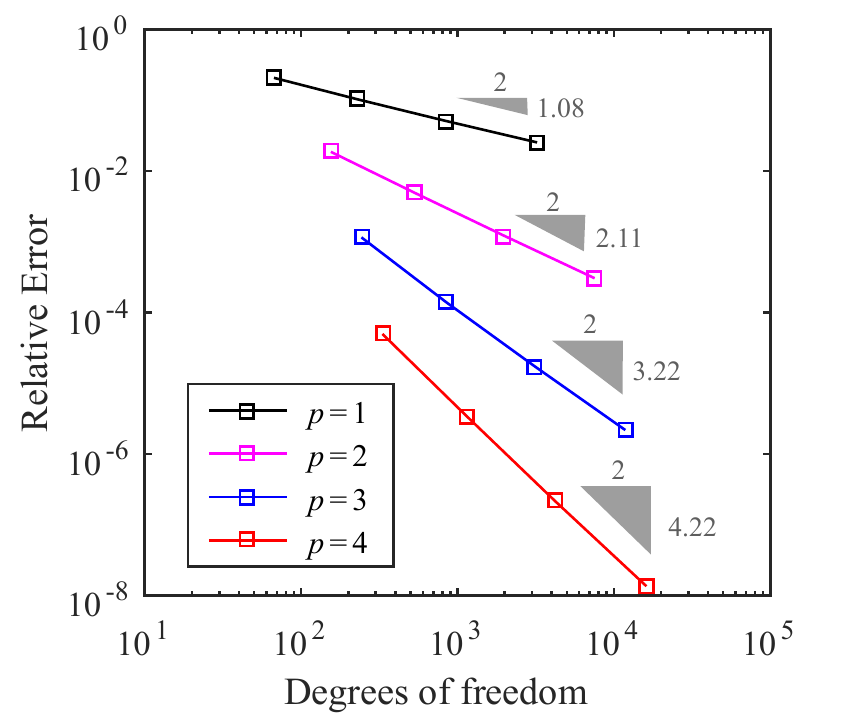}
	\caption{Convergence study of the PolyDPG method in terms of degrees of freedom. The $h$-convergence behavior is displayed for different $p$ using the meshes with an interface in Figure \ref{fig:fault_mesh}.}
	\label{fig:fault_error}
\end{figure}

Figure \ref{fig:fault_plots} shows the appearance of the computed ultraweak solution. 
As it can be observed in Figure~\ref{fig:fault_error}, the expected convergence rates were verified once again. 
It is remarkable that without collapsing any nodes in these meshes, the same data points were observed for $1\leq p\leq3$, but the last data point for $p=4$ did behave unexpectedly, so collapsing the nodes is still recommended in general. 

\subsection{Adaptivity} 
\label{sec:results_adapt}

In the last example, we aimed to present a polygonal adaptive strategy.
This is of interest as it has direct applications in fracture dynamics \cite{pfem_cracks,fracturemeshbias} and topology optimization \cite{pfem_optimization,vem_optimization}.
Implementing such a strategy was possible, because the DPG methodology carries a natural arbitrary-order \textit{a posteriori} error estimator (see \eqref{eq:residualexpression} and Section \ref{sec:DPGcharacteritics}).
Indeed, assuming that $\eta_K$ is the \textit{a posteriori} error estimator (representing the square of the residual as in \eqref{eq:residualexpression}) for $K\in\mesh$, and $\eta_{\max}=\max_{K\in\mesh}\eta_K$, then the criterion used to mark an element for refinement was if $\eta_K\geq0.25\eta_{\max}$ \cite{demkowicz2012class}.

In order to refine traditional quadrilateral elements, typically hanging nodes arise in the mesh. 
But in practice, only one ``level'' of refinement is possible per element (often edges cannot have more than one hanging node), resulting in so-called quadtree meshes \cite{pfem_quadtree}.
To implement this strategy a constrained approximation technology is necessary to handle the hanging nodes.
Additionally, under anisotropic refinements, sometimes dead-lock scenarios arise (where it is logically impossible to continue refining) and these must be avoided \cite{hpbook}.
In short, it may be challenging to implement conventional refinement strategies used for adaptivity.

\begin{figure}[!ht]
	\centering
	\includegraphics[width=13cm]{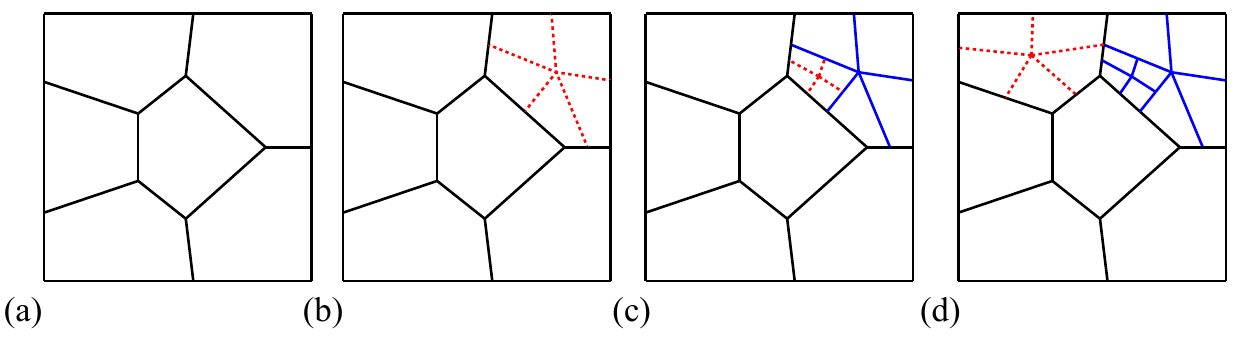}
	\caption{A practical local refinement strategy for convex polygons: (a) initial coarse polygonal mesh; (b) line segments are projected from the centroid to every edge midpoint in the element of interest; (c) the same approach is used to refine sub-elements; (d) the strategy can be re-applied to any other coarser element by assuming all collinear vertices constitute a single combined edge.}
	\label{fig:local_refine}
\end{figure}

An important advantage of the polygonal elements is that they naturally embrace hanging nodes, because they merely represent that a polygon has an extra edge collinear with another edge.
Thus, the polygonal methods do not require an extra level of difficulty in terms of implementing the adaptive refinements.
We devised a practical convex polygonal refinement strategy as illustrated in Figure \ref{fig:local_refine}: (a) shows the initial mesh in which an element of interest is picked and split into quadrilaterals by using the centroid and edge midpoints as depicted in (b); next, any of the resulting elements can be subsequently refined into finer quadrilaterals as shown in (c); and lastly, as shown in (d), if a neighbor element needs to be refined too, it is split into quadrilaterals assuming all adjacent collinear edges constitute a single edge (i.e.~the vertices of this combined edge are used in the calculation of the centroid and its midpoint used to place the new quadrilateral node).

%
%

The manufactured solution for this problem is the sum of two Gaussian surfaces, given by the function,
\begin{equation}
	\label{eq:manuf4}
	\centering
	u(x,y) = \frac{1}{2\pi\sigma^2}\left[e^{-\frac{1}{2}\left(\frac{x-\mu_1}{\sigma}\right)^2}e^{-\frac{1}{2}
		\left(\frac{y-\mu_1}{\sigma}\right)^2}+e^{-\frac{1}{2}
			\left(\frac{x-\mu_2}{\sigma}\right)^2}e^{-\frac{1}{2}\left(\frac{y-\mu_2}{\sigma}\right)^2}\right]\ ,
\end{equation}
where the standard deviation is $\sigma=\sqrt{10^{-3}}$ and the two means are $\mu_1=0.25$ and $\mu_2=0.75$. 
Again, projection-based interpolation \cite{demkowicz2008polynomial,hpbook2} was used to approximate the nearly vanishing temperature boundary conditions.

\begin{figure}[!ht]
	\centering
	\includegraphics[width=16cm]{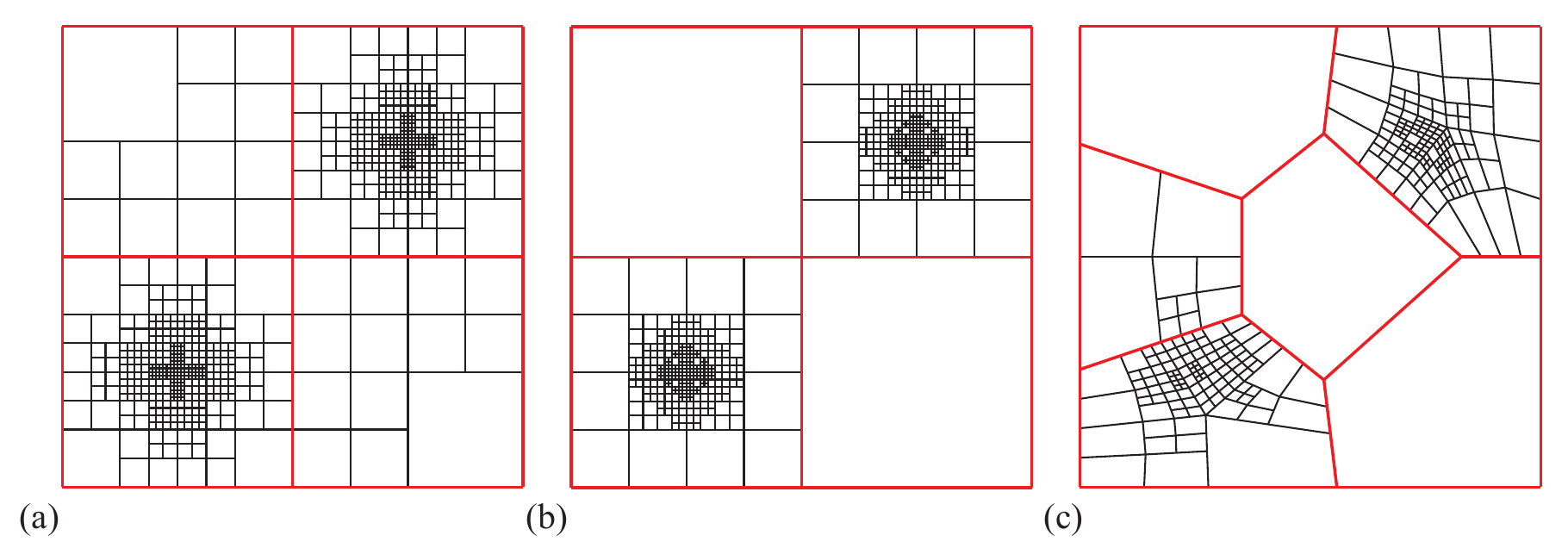}
	\caption{Three $h$-adaptively refined meshes (the red line represents the initial mesh) for the manufactured solution in \eqref{eq:manuf4}: (a) traditional quadtree meshes via constrained nodes; (b) quadrilateral mesh using the polygonal adaptive strategy; (c) polygonal mesh using the polygonal adaptive strategy.}
	\label{fig:adapt_mesh}
\end{figure}

In order to compare with other adaptive schemes, a traditional adaptive strategy using quadtree meshes and constrained hanging nodes via quadrilateral elements was considered here \cite{hpbook}.
Starting with the same initial mesh, the traditional refinement strategy and the polygonal refinement strategy were allowed to refine accordingly.
When using the polygonal strategy on these quadrilateral meshes, we used the more natural choice of bounding boxes instead of the bounding triangles.
Additionally, the same polygonal refinement strategy was applied to an initial polygonal mesh (using bounding triangles as usual).
Figure \ref{fig:adapt_mesh} shows the results of the three different scenarios after several refinements.
Clearly, the traditional adaptive strategy produces quadtree meshes (see Figure \ref{fig:adapt_mesh}(a)), so it is forced to refine and create new elements in areas of the domain where the solution is nearly constant.
However, the polygonal adaptive strategy applied to the same initial mesh produces a more localized refinement pattern which is not a quadtree mesh (see Figure \ref{fig:adapt_mesh}(b)).
Lastly, the polygonal adaptive strategy applied to a polygonal mesh produces a completely nonstandard, yet localized mesh (see Figure \ref{fig:adapt_mesh}(c)).

\begin{figure}[!ht]
	\centering
	\includegraphics[width=16cm]{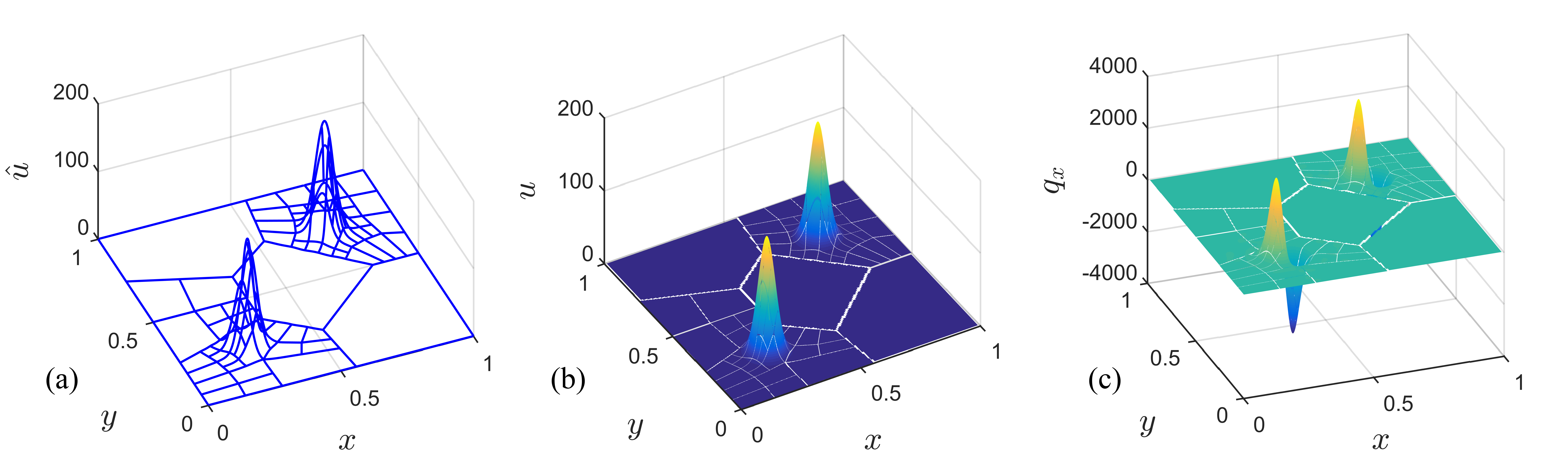}
	\caption{Numerical results using the manufactured solution in \eqref{eq:manuf4} on the mesh from Figure \ref{fig:adapt_mesh}(c) using $p=6$ and $\Delta p=2$: (a) skeleton temperature, (b) temperature, (c) first component of the heat flux.}
	\label{fig:adapt_plots}
\end{figure}

\begin{figure}[!ht]
	\centering
	\includegraphics[width=8cm]{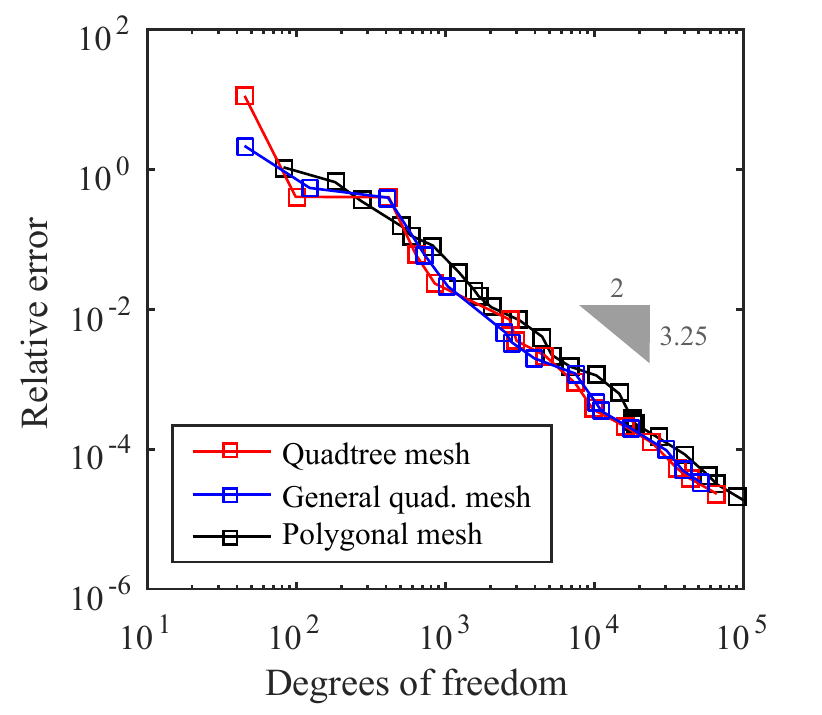}
	\caption{Convergence study of the PolyDPG method in terms of degrees of freedom. The $h$-convergence behavior is displayed using $p=3$ for several successive refinements associated with the refinement strategies in Figure \ref{fig:adapt_mesh}.}
	\label{fig:adapt_error}
\end{figure}

The numerical solution for $p=6$ and $\Delta p=2$ using the mesh in Figure \ref{fig:adapt_mesh}(c) is presented in Figure~\ref{fig:adapt_plots}.
The error convergence curves corresponding to the three refinement schemes in Figure \ref{fig:adapt_mesh} are also displayed in Figure~\ref{fig:adapt_error}.
The proposed polygonal refinement technique generates more edges (each new sub-segment becomes an edge) resulting in more degrees of freedom. However, in the end the additional cost is compensated by producing less elements than traditional quadtree refinement schemes (compare (b) and (c) with (a) in Figure \ref{fig:adapt_mesh}).
It can be seen from Figure \ref{fig:adapt_error} that the convergence behavior in terms of degrees of freedom is very similar using both approaches.
Therefore, the polygonal adaptive strategy proposed here is competitive with the existing strategies for traditional elements, whilst being more general in its applicability as it also works for polygonal elements. 

%% file: 4_Conclusions.tex
\section{Conclusions} 
\label{sec:conclusions}

A PolyDPG method discretized with high-order polygonal elements was successfully implemented using ultraweak formulations and the DPG methodology. 
Here, the PolyDPG method solves Poisson's equation. 
However, like with the DPG methodology, the discretization and theory is quite general.
Thus, it can be applied to a large family of equations including acoustics, convection-dominated diffusion and linear elasticity.
PolyDPG methods are conforming FEMs, and as with many other polytopal methods, the spaces and integration schemes are defined directly in the physical space.
Indeed, given that the ultraweak formulations avoid inter-element compatibility conditions, it is relatively straightforward to obtain many of the shape functions by restricting them from a bounding (triangular or quadrilateral) element to the polygonal element.
Despite the greater computational cost compared to conventional methods, the resulting PolyDPG methods are naturally high-order, carry their own residual-based \textit{a posteriori} error estimator, have no need of ad hoc stabilization terms, and always produce positive-definite stiffness matrices. 
Moreover, under reasonable assumptions, a rigorous proof demonstrating the convergence of PolyDPG methods was included.
To complement this work, the \texttt{PolyDPG} software \cite{polydpg} written in MATLAB\textsuperscript{\textregistered} is provided.
We hope this will prove to be a practical tool for other researchers interested in polygonal FEMs and in DPG methods. 

Different illustrative examples corroborated the expected results. 
In the first example, $n$-sided convex polygons were investigated, while in the second example, highly distorted concave elements were examined.
In both cases, as predicted by the theory, convergence rates of the form $h^p$ were observed for different values of $p$, confirming that PolyDPG methods are distortion-tolerant.
The third example was relevant to the field of geosciences, where faults cause heterogeneity in the domain.
This was simulated by irregularly cutting a uniform grid with an interface and assigning different material properties on each side.
Once again, the method converged as expected, displaying its robustness in resolving heterogeneous material properties.  
The final example explored a polygonal adaptivity scheme driven by the arbitrary-order \textit{a posteriori} error estimator of PolyDPG methods.
Even though polygonal and standard refinement strategies led to practically identical convergence curves, polygonal techniques are more general since they apply to polygonal elements and avoid the typical approaches of constrained approximations via hanging nodes.  
These techniques may be useful in applications such as crack propagation and topology optimization.

Extension of the presented technique to arbitrary 3D polyhedral elements is in progress. 
In principle, the current numerical method can be extended naturally to polyhedral elements, as long as all the faces are triangular, but the case of arbitrary faces is much more challenging and might lead to analyzing nonconforming numerical methods.

\paragraph{Acknowledgments.}
This work was partially supported with grants by ONR (N00014-15-1-2496), NSF (DMS-1418822), and AFOSR (FA9550-12-1-0484).
Co-author Jaime Mora is also sponsored by a 2015 Colciencias-Fulbright scholarship (Colombia).

%% file: 5_Convergence.tex
\section{Convergence} 
\label{sec:AppConvergence}

\subsection{Stability and Fortin operators}

Since the numerical method is technically a conforming FEM, the ``standard'' theory of convergence can be applied.
However, the issue of numerical stability, in the sense of \eqref{eq:infsupstable}, must be addressed first.
The DPG methodology is basically crafted to ``almost'' satisfy this condition, and intuitively, the larger the enriched test space $\scV_r\subseteq\scV$, the more certainty there is that the condition is satisfied.
This translates to increasing $\Delta p_K$ for all $K\in\mesh$ in \eqref{eq:Vrdefinition}, so that $\scV_r$ becomes larger. Note that this increases the local (element-wise) computational burden.
In practice, the numerical stability is observed even with very modest values of $\Delta p$.

However, to have a rigorous result, it is necessary to establish \eqref{eq:infsupstable} theoretically.
To do so, it is helpful to consider a linear and continuous Fortin operator, $\Pi_F:\scV\to\scV_r$, satisfying the orthogonality condition, $b(\fku_h,\fkv-\Pi_F\fkv)=0$, for all $\fku_h\in\scU_h$ and $\fkv\in\scV$.
If it exists, it follows that \cite{QiuFortin},
\begin{equation}
	\|\fku-\fku_h\|_\scU\leq\frac{\|b\|\,M_F}{\gamma}\inf_{\delta\fku_h\in U_h}\|\fku-\delta\fku_h\|_\scU\,,
	\label{eq:convergencewithFortin}
\end{equation}
where $M_F\geq\|\Pi_F\|=\sup_{\fkv\in \scV}\frac{\|\Pi_F\fkv\|_\scV}{\|\fkv\|_\scV}$, $\|b\|=\sup_{(\fku,\fkv)\in \scU\!\times\! \scV}\frac{|b(\fku,\fkv)|}{\|\fku\|_\scU\|\fkv\|_\scV}$ and $\gamma=\inf_{\fku\in\scU}\sup_{\fkv\in\scV}\frac{|b(\fku,\fkv)|}{\|\fku\|_\scU\|\fkv\|_\scV}$, where the infima and suprema are tacitly assumed to be taken over nonzero elements.
Note that when $\scV$ is a broken test space, as in this case, the Fortin operator can be separately constructed locally at each element $K\in\mesh$.
Constructions of such Fortin operators do exist for triangles \cite{BrokenForms15}, but have not been constructed for other shapes yet.
Nevertheless, numerical results show it is reasonable to expect them to exist, and this will be assumed in what follows.
In any case, note that Fortin operators merely yield a conservative estimate, but in practice the results are better (i.e. instead of $M_F$, there is a moderate constant, $\mcO(1)$-$\mcO(10)$, multiplying $\frac{\|b\|}{\gamma}$ in \eqref{eq:convergencewithFortin}).

\subsection{Fractional spaces}

For a given $s\geq0$ and any Lipschitz domain $A\subseteq\R^2$, the fractional Sobolev spaces, $\HSo^{1+s}(A)$, $\bHSo^s(\div,A)$ and $\HSo^s(A)$, are slightly smoother subspaces of $\HSo^1(A)$, $\bHSo(\div,A)=\bHSo^0(\div,A)$ and $\Leb^2(A)=\HSo^0(A)$, respectively.
As an obvious placeholder, let $\scW(A)=\scW^0(A)$ be one of these spaces, and for $s\geq0$, let $\scW^s(A)$ be its slightly smoother fractional counterpart.
As one might expect, $\|w\|_{\scW(A)}\leq\|w\|_{\scW^{s_1}(A)}\leq\|w\|_{\scW^{s_2}(A)}$ for $0\leq s_1\leq s_2$ and all $w\in\scW^{s_2}(A)\subseteq\scW^{s_1}(A)\subseteq\scW(A)$.
For more details on these spaces and their norms, see \cite{McLeanSobolev}.

Using interpolation theory (see \cite[Appendix B]{McLeanSobolev}) applied to the universal extension operators of Sobolev spaces of differential forms defined in \cite{UniversalExtensionSobolev} (which is even more general than the universal extension operator defined by Stein in \cite{SteinUniiversalExtension}), it is possible to establish the existence of a continuous extension operator,
\begin{equation}
	E:\scW^s(\Omega)\to\scW^s(\R^2)\,,\qquad \|Ew\|_{\scW^s(\R^2)}\leq C_E\|w\|_{\scW^s(\Omega)}\,,
	\label{eq:extensionoperator}
\end{equation}
where $s\geq0$, $\Omega$ is the domain where the equations are being solved, and $C_E=C_E(s,\Omega)>0$. 

The fractional skeleton spaces are better defined directly through fractional trace operators of Lipschitz elements $K\in\mesh$ (see \cite{McLeanSobolev}) as $\HSo^{\onehalf+s}(\bdry K)=\tr_{\HSo^{1+s}(K)}(\HSo^{1+s}(K))$ and $\HSo^{-\onehalf+s}(\bdry K)=\tr_{\bHSo^s(\div,K)}(\bHSo^s(\div,K))$ (see \eqref{eq:localtracediscretizations} for the explicit trace operators for the $s=0$ case).
Again, using placeholders these are written as $\scW^s(\bdry K)=\{\hat{w}_K=\tr_{\scW^s(K)}w\mid w\in\scW^s(K)\}$, so that their minimum energy extension norm is
\begin{equation}
	\|\hat{w}_K\|_{\scW^s(\bdry K)}=\inf_{w\in\tr_{\scW^s(K)}^{-1}\{\hat{w}_K\}}\|w\|_{\scW^s(K)}\,.
	\label{eq:minenergyboundary}
\end{equation}
At the global level, define the global trace operators as
\begin{equation}
	\tr_{\scW^s(\mesh)}:\scW^s(\Omega)\to\prod_{K\in\mesh}\scW^s(\bdry K)\,,\qquad
		\tr_{\scW^s(\mesh)}w=\prod_{K\in\mesh}\tr_{\scW^s(K)}(w|_K)\,.
\end{equation}
Note that $\HSo_0^{1+s}(\Omega)=\overbar{C_0^\infty(\Omega)}{}^{\|\cdot\|_{\HSo^{1+s}(\Omega)}}$, so that the global fractional skeleton spaces are (see \eqref{eq:interfacespaces} for the $s=0$ case),
\begin{equation}
	\HSo_0^{\onehalf+s}(\bdry\mesh)=\tr_{\HSo^{1+s}(\mesh)}\big(\HSo_0^{1+s}(\Omega)\big)\,,\qquad
		\HSo^{-\onehalf+s}(\bdry\mesh)=\tr_{\bHSo^s(\div,\mesh)}\big(\bHSo^s(\div,\Omega)\big)\,.
\end{equation}

Analogous to \eqref{eq:BrokenUWSpaces}, the fractional trial subspace for $s\geq0$ is 
\begin{equation}
	\scU^s=\HSo^s(\Omega)\times(\HSo^s(\Omega))^2\times
		\HSo_0^{\onehalf+s}(\bdry\mesh)\times\HSo^{-\onehalf+s}(\bdry\mesh)\subseteq\scU\,,
\end{equation}
and it is easy to see $\|\fku\|_{\scU}\leq\|\fku\|_{\scU^{s_1}}\leq\|\fku\|_{\scU^{s_2}}$ for $0\leq s_1\leq s_2$ and all $\fku\in\scU^{s_2}\subseteq\scU^{s_1}\subseteq\scU$.

\subsection{Approximation properties}

Next, for a bounded Lipschitz domain $A\subseteq\R^2$ and polynomial order $p\in\N$, consider commuting exact sequence discretizations for $\HSo^1(A)$, $\bHSo(\div,A)$ and $\Leb^2(A)$, such that
\begin{equation}
	\begin{gathered}
		\xymatrix@R-1.5pc{
		\HSo^1(A) &\bHSo(\div,A) &\Leb^2(A)\\
		\HSo_{hp}^1(A) \ar@{}[u]|*=0[@]{\subseteq} \ar[r]^{\curl\quad\,\,} 
			&\bHSo_{hp}(\div,A) \ar@{}[u]|*=0[@]{\subseteq} \ar[r]^{\quad\nabla\cdot} 
				&\Leb_{hp}^2(A) \ar@{}[u]|*=0[@]{\subseteq}\\
		\mathcal{P}^p(A) \ar@{}[u]|*=0[@]{\subseteq}
			&(\mathcal{P}^{p-1}(A))^2 \ar@{}[u]|*=0[@]{\subseteq} 
				&\mathcal{P}^{p-1}(A) \ar@{}[u]|*=0[@]{\subseteq}\,.
		}
	\end{gathered}
	\label{eq:exactsequenceconvergence}
\end{equation}
More abstractly, the discretizations are written as $\scW_{hp}(A)\subseteq\scW(A)$.
Then, given the polynomials contained in the discretizations $\scW_{hp}(A)$, it is well known that for $p\geq s>0$, there exists a constant $C_h=C_h(A,s)>0$, such that for all $w\in\scW^s(A)$,
\begin{equation}
	\inf_{\delta w\in\scW_{hp}(A)}\|w-\delta w\|_{\scW^s(A)}\leq C_h\|w\|_{\scW^s(A)}\,.
	\label{eq:simpleapprox}
\end{equation}
For each $K\in\mesh$, the local trace spaces are supposed to be $\scW_{hp}(\bdry K)=\tr_{\scW^0(K)}(\scW_{hp}(K))$ for some $\scW_{hp}(K)$.

We would like these approximation properties to hold for our choices of discrete trial spaces, $\scU_h$ in \eqref{eq:Uhdefinition}, and this is indeed the case.
The first two components of $\scU_h$ when restricted to $K\in\mesh$ (representing $\Leb_{hp}^2(K)$ and $(\Leb_{hp}^2(K))^2$) are restrictions to $K$ of $\mcP^{p-1}(T_K)$ and $(\mcP^{p-1}(T_K))^2$, so they do trivially contain $\mcP^{p-1}(K)$ and $(\mcP^{p-1}(K))^2$ respectively.
This means that \eqref{eq:simpleapprox} holds for those two spaces, but as we will see soon, it suffices (and is preferable) to have this result for the bounding triangle $T_K$ (which is obviously true).
For the third and fourth components of $\scU_h$, representing the skeleton variables, locally at each  $K\in\mesh$ it suffices to show that $\mcP^{p}_C(\bdry K)=\tr_{\HSo^1(K)}(\HSo_{hp}^1(K))$ and $\mcP^{p-1}(\bdry K)=\tr_{\bHSo(\div,K)}(\bHSo_{hp}(\div,K))$ for some $\HSo_{hp}^1(K)$ and $\bHSo_{hp}(\div,K)$ satisfying the properties in \eqref{eq:exactsequenceconvergence}, where $\mcP^{p}_C(\bdry K)$ and $\mcP^{p-1}(\bdry K)$ are defined in \eqref{eq:localtracediscretizations}.
For this, consider the shape-regular edge-compatible triangulations of each $K\in\mesh$, denoted by $\scT(K)=\{\scT_i(K)\}_{i\in I_K}$ (with $I_K$ finite), and define the spaces,
\begin{equation}
	\begin{aligned}
		\HSo_{hp}^1(K)&=\{u\in\HSo^1(K)\mid u|_{\scT_i(K)}\in\mcP^p(\scT_i(K)),\forall i\in I_K\}\,,\\
		\bHSo_{hp}(\div,K)&=\{\bsq\in\bHSo(\div,K)\mid\bsq|_{\scT_i(K)}\in\mcR\mcT^p(\scT_i(K)),\forall i\in I_K\}\,.
	\end{aligned}
	\label{eq:discretetracesubmesh}
\end{equation}
It can easily be checked that $\hat{u}_K=u|_{\bdry K}\in\mcP^{p}_C(\bdry K)$ for all $u\in\HSo_{hp}^1(K)$ and $(\hat{q}_{\nml})_K=\bsq|_{\bdry K}\!\cdot\!\nml_K\in\mcP^{p-1}(\bdry K)$ for all $\bsq\in\bHSo_{hp}(\div,K)$, and that these inclusions are surjective.
Thus, $\mcP^{p}_C(\bdry K)=\tr_{\HSo^1(K)}(\HSo_{hp}^1(K))$ and $\mcP^{p-1}(\bdry K)=\tr_{\bHSo(\div,K)}(\bHSo_{hp}(\div,K))$ as desired.
This implies that \eqref{eq:simpleapprox} also holds for $\HSo_{hp}^1(K)$ and $\bHSo_{hp}(\div,K)$, which are closely related to the skeleton discretizations of $\scU_h$. 

\subsection{Interpolation estimates}

The idea is to define a bounded linear interpolation operator $\Pi_{\scU^s}:\scU^s\to\scU_h$ such that $\Pi_{\scU^s}\fku_h=\fku_h$ for every $\fku_h\in\scU_h$ and $s>\frac{1}{2}$.
Typically this implies constructing interpolation operators for every component of $\scU$. Moreover, for each component this construction is done locally at every $K\in\mesh$ in such a way that the inter-element compatibility properties are satisfied.

The first two components of $\scU$ are $\Leb^2(\Omega)$ and $\bLeb^2(\Omega)$, which are effectively three $\Leb^2(\Omega)$ components. 
The last two skeleton components are $\HSo_0^{\onehalf+s}(\bdry\mesh)$ and $\HSo^{-\onehalf+s}(\bdry\mesh)$.
The discretizations of these three spaces are (see \eqref{eq:Uhdefinition}),
\begin{equation}
	\begin{aligned}
		\Leb_{hp}^2(\Omega)&=\big\{u\in\Leb^2(\Omega)\mid u|_K\in\mcP^{p-1}(K),\forall K\in\mesh\big\}\,,\\
		\HSo_{0,hp}^{\onehalf}(\bdry\mesh)&=\big\{\hat{u}\in\HSo_0^{\onehalf}(\bdry\mesh)
			\mid\hat{u}_K\in\mcP^{p}_C(\bdry K)=\tr_{\HSo^1(K)}(\HSo_{hp}^1(K)),\forall K\in\mesh\big\}\,,\\
		\HSo_{hp}^{-\onehalf}(\bdry\mesh)&=\big\{\hat{q}_{\nml}\in\HSo^{-\onehalf}(\bdry\mesh)
			\mid(\hat{q}_{\nml})_K\in\mcP^{p-1}(\bdry K)=\tr_{\bHSo(\div,K)}(\bHSo_{hp}(\div,K)),\forall K\in\mesh\big\}\,,
	\end{aligned}
\end{equation}
where the definitions of $\HSo_{hp}^1(K)$ and $\bHSo_{hp}(\div,K)$ are in \eqref{eq:discretetracesubmesh}.
Thus, it suffices to construct,
\begin{equation}
	\begin{aligned}
		\Pi_{\HSo^s(\Omega)}:\HSo^s(\Omega)\to\Leb_{hp}^2(\Omega)\,,&\qquad
			\big(\Pi_{\HSo^s(\Omega)}u\big)\big|_K=\Pi_{\HSo^s(K)}u|_K\,,\\
		\Pi_{\HSo_0^{\onehalf+s}(\bdry\mesh)}:\HSo_0^{\onehalf+s}(\bdry\mesh)\to\HSo_{0,hp}^{\onehalf}(\bdry\mesh)\,,&\qquad
			\big(\Pi_{\HSo_0^{\onehalf+s}(\bdry\mesh)}\hat{u}\big)_K=\Pi_{\HSo^{\onehalf+s}(\bdry K)}\hat{u}_K\,,\\
		\Pi_{\HSo^{-\onehalf+s}(\bdry\mesh)}:\HSo^{-\onehalf+s}(\bdry\mesh)\to\HSo_{hp}^{-\onehalf}(\bdry\mesh)\,,&\qquad
			\big(\Pi_{\HSo^{-\onehalf+s}(\bdry\mesh)}\hat{q}_{\nml}\big)_K=\Pi_{\HSo^{-\onehalf+s}(\bdry K)}(\hat{q}_{\nml})_K\,,
	\end{aligned}
	\label{eq:globalinterpolants}
\end{equation}
meaning that we must define $\Pi_{\HSo^s(K)}$, $\Pi_{\HSo^{\onehalf+s}(\bdry K)}$ and $\Pi_{\HSo^{-\onehalf+s}(\bdry K)}$.

The operator $\Pi_{\HSo^s(K)}$ can be chosen as the $\Leb^2(K)$-projection to $\mcP^{p-1}(K)$ directly on $K$ (so $\Pi_{\HSo^s(K)}\delta u=\delta u$ for all $\delta u\in\mcP^{p-1}(K)$).
Consider now a simple scaling by $h_K=\diam(K)$, so that $\hat{K}$ has $\diam(\hat{K})=1$.
Using \eqref{eq:simpleapprox} for $p\geq s>\frac{1}{2}$ results in the abstract expression,
\begin{equation}
	\begin{aligned}
		\|w-\Pi_{\scW^s(\hat{K})}w\|_{\scW(\hat{K})}
			&=\inf_{\delta w\in\scW_{hp}(\hat{K})}\|(I-\Pi_{\scW^s(\hat{K})})(w-\delta w)\|_{\scW(\hat{K})}\\
				&\leq\|I-\Pi_{\scW^s(\hat{K})}\|\inf_{\delta w\in\scW_{hp}(\hat{K})}\|(w-\delta w)\|_{\scW^s(\hat{K})}
					\leq C_{\hat{K}}\|w\|_{\scW^s(\hat{K})}\,,
	\end{aligned}
\end{equation}
for any $w\in\scW^s(\hat{K})$, where $C_{\hat{K}}=C_{\hat{K}}(\hat{K},p,s)>0$.
Scaling appropriately then yields for any $w\in\scW^s(K)$,
\begin{equation}
	\|w-\Pi_{\scW^s(K)}w\|_{\scW(K)}\leq C_{\hat{K}}h_K^s\|w\|_{\scW^s(K)}\,.
	\label{eq:hsestimatelocaldirect}
\end{equation}
The issue with this estimate is that it depends on the element shape $K$ (via $\hat{K}$), so it is inconvenient as it may become much larger with mesh refinements.
The solution is to use the bounding triangle and the extension operator defined in \eqref{eq:extensionoperator}, so that (as in \cite{hpdgfem}) the interpolation operator is defined for any $w\in\scW^s(\Omega)$ as,
\begin{equation}
	\Pi_{\scW^s(K)}w|_K=\big(\Pi_{\scW^s(T_K)}Ew|_{T_K}\big)\big|_K\,,
\end{equation}
where $\Pi_{\scW^s(T_K)}$ is the $\Leb^2(T_K)$-projection.
Scaling and rotating transforms the bounding triangle $T_K$ to a unique triangle $\hat{T}_0$ (independent of the element $K$) with $\diam(\hat{T}_0)=1$.
This means $T_K$ is scaled by $h_{T_K}=\diam(T_K)=\frac{6}{\sqrt{3}}r_{\max}\leq\sqrt{12}h_K$, where $r_{\max}$ is the distance of the centroid to the furthest vertex and $h_K=\diam(K)$.
Using the same reasoning gives,
\begin{equation}
	\|w-\Pi_{\scW^s(K)}w\|_{\scW(K)}\leq\|Ew-\Pi_{\scW^s(T_K)}Ew\|_{\scW(T_K)}\leq C_{\hat{T}_0}h_K^s\|Ew\|_{\scW^s(T_K)}\,,
\end{equation}
for every $w\in\scW^s(K)$, where $C_{\hat{T}_0}=C_{\hat{T}_0}(p,s)>0$ is now independent of $K$.

Next, consider the skeleton variables for an element $K\in\mesh$ and its respective shape-regular and edge-compatible triangulation denoted by $\scT(K)=\{\scT_i(K)\}_{i\in I_K}$.
The theory of projection-based interpolation \cite{demkowicz2008polynomial} implies that for any polygonal domain $A\subseteq\R^2$ and for any $s>\frac{1}{2}$ there exist commuting operators,
\begin{equation}
	\begin{gathered}
  	\xymatrix{
        & \HSo^{1+s}(A) \ar[r]^{\curl\quad} \ar@{>}[d]^{\Pi_{\HSo^{1+s}(A)}}
          & \bHSo^s(\div,A) \ar[r]^{\,\,\,\nabla\cdot} \ar@{>}[d]^{\Pi_{\bHSo^s(\div,A)}} 
						& \HSo^s(A) \ar@{>}[d]^{\Pi_{\HSo^s(A)}}\\
        & \HSo_{hp}^1(A) \ar[r]^{\curl\quad} 
					& \bHSo_{hp}(\div,A)\ar[r]^{\,\,\,\nabla\cdot} 
						& \Leb_{hp}^2(A)\,.
						}
	\end{gathered}
	\label{eq:projectionbasedsequences}
\end{equation}
Thus, for any $p\geq s>\frac{1}{2}$ and triangle $\scT_i(K)$ (so $\diam(\scT_i(K))\leq\diam(K)=h_K$), the result in \eqref{eq:hsestimatelocaldirect} applies and yields
\begin{equation}
	\|w-\Pi_{\scW^s(\scT_i(K))}w\|_{\scW(\scT_i(K))}\leq C_{\hat{\scT}_0}h_K^s\|w\|_{\scW^s(\scT_i(K))}\,,
\end{equation}
where the $K$-independent $C_{\hat{\scT}_0}=C_{\hat{\scT}_0}(p,s)>0$ exists due to the \textit{assumed} uniform shape-regularity of the $\scT_i(K)$ (across all $K\in\mesh$ and all meshes being considered).
Adding among $\scT(K)$ is valid due to the compatibility of the projection-based interpolation in the triangulation, so that
\begin{equation}
	\|w-\Pi_{\scW^s(K)}w\|_{\scW(K)}^2=\sum_{i\in I_K}\|w-\Pi_{\scW^s(\scT_i(K))}w\|_{\scW(\scT_i(K))}^2
		\leq C_{\hat{\scT}_0}^2h_K^{2s}\|w\|_{\scW^s(K)}^2\,.
\end{equation}
Lastly, consider the well-defined trace interpolation,
\begin{equation}
	\Pi_{\scW^s(\bdry K)}\hat{w}_K=\tr_{\scW(K)}\Pi_{\scW^s(K)}w\,,\qquad w\in\tr_{\scW^s(K)}^{-1}\{\hat{w}_K\}\,,
\end{equation}
so that (see \eqref{eq:minenergyboundary}),
\begin{equation}
	\begin{aligned}
		\|\hat{w}_K-\Pi_{\scW^s(\bdry K)}\hat{w}_K\|_{\scW(\bdry K)}&=\|\tr_{\scW(K)}(w-\Pi_{\scW^s(K)}w)\|_{\scW(\bdry K)}\\
			&\leq\|w-\Pi_{\scW^s(K)}w\|_{\scW(K)}\leq C_{\hat{\scT}_0}h_K^s\|w\|_{\scW^s(K)}\,.
	\end{aligned}
\end{equation}
This is true for every $w\in\tr_{\scW^s(K)}^{-1}\{\hat{w}_K\}$, so take the infimum to yield
\begin{equation}
		\|\hat{w}_K-\Pi_{\scW^s(\bdry K)}\hat{w}_K\|_{\scW(\bdry K)}\leq C_{\hat{\scT}_0}h_K^s\|\hat{w}_K\|_{\scW^s(\bdry K)}\,.
\end{equation}

Putting everything together and generalizing for any $p\in\N$ and $s>\frac{1}{2}$, gives
\begin{equation}
	\begin{aligned}
		\|u-\Pi_{\HSo^s(K)}u\|_{\Leb^2(K)}&\leq C_{\hat{T}_0}h_K^{\min\{s,p\}}\|Eu\|_{\HSo^s(T_K)}\,,\\
		\|\bsq-\Pi_{(\HSo^s(K))^2}\bsq\|_{\bLeb^2(K)}&\leq C_{\hat{T}_0}h_K^{\min\{s,p\}}\|E\bsq\|_{(\HSo^s(T_K))^2}\,,\\
		\|\hat{u}_K-\Pi_{\HSo^{\onehalf+s}(\bdry K)}\hat{u}_K\|_{\HSo^{\onehalf}(\bdry K)}
			&\leq C_{\HSo^{1+s}(\hat{\scT}_0)}h_K^{\min\{s,p\}}\|\hat{u}_K\|_{\HSo^{\onehalf+s}(\bdry K)}\,,\\
		\|(\hat{q}_{\nml})_K-\Pi_{\HSo^{-\onehalf+s}(\bdry K)}(\hat{q}_{\nml})_K\|_{\HSo^{-\onehalf}(\bdry K)}
			&\leq C_{\bHSo^s(\div,\hat{\scT}_0)}h_K^{\min\{s,p\}}\|(\hat{q}_{\nml})_K\|_{\HSo^{-\onehalf+s}(\bdry K)}\,,		
	\end{aligned}
	\label{eq:hestimatesconvergence}
\end{equation}
where the constants $C_{\hat{T}_0}$, $C_{\HSo^{1+s}(\hat{\scT}_0)}$ and $C_{\bHSo^s(\div,\hat{\scT}_0)}$ only depend on $p$ and $s$, but not on $K$ (the last two constants depend on the uniform shape-regularity of the edge-compatible triangulations of all elements).
Finally, since these constants come from triangles, the theory of projection-based interpolation \cite{demkowicz2008polynomial} implies that in the $p$-asymptotic limit,
\begin{equation}
	C_{\hat{T}_0}=\tilde{C}_{\hat{T}_0}(\ln p)p^{-s}\,,\quad\,\, 
		C_{\HSo^{1+s}(\hat{\scT}_0)}=\tilde{C}_{\HSo^{1+s}(\hat{\scT}_0)}(\ln p)^2p^{-s}\,,\quad\,\,
			C_{\bHSo^s(\div,\hat{\scT}_0)}=\tilde{C}_{\bHSo^s(\div,\hat{\scT}_0)}(\ln p)p^{-s}\,,
	\label{eq:pdependenceofC}
\end{equation}
where $\tilde{C}_{\hat{T}_0}$, $\tilde{C}_{\HSo^{1+s}(\hat{\scT}_0)}$ and $\tilde{C}_{\bHSo^s(\div,\hat{\scT}_0)}$ are constants independent of $p$ and of any $K\in\mesh$ across all possible meshes being considered.

\subsection{Final convergence estimates}

Use the global interpolation operators in \eqref{eq:globalinterpolants} to construct the bounded linear global interpolation operator $\Pi_{\scU^s}:\scU^s\to\scU_h$.
Note that adding \eqref{eq:hestimatesconvergence} associated with $u\in\HSo^s(\Omega)$ among $K\in\mesh$, using the robust finite overlap condition, and the extension operator in \eqref{eq:extensionoperator}, gives:
\begin{equation}
	\begin{aligned}
		\|u-\Pi_{\HSo^s(\Omega)}u\|_{\Leb^2(\Omega)}^2&\leq C_{\hat{T}_0}^2\sum_{K\in\mesh}h_K^{2\min\{s,p\}}\|Eu\|_{\HSo^s(T_K)}^2\\
			&\leq M_{\overlap}C_{\hat{T}_0}^2h^{2\min\{s,p\}}\|Eu\|_{\HSo^s(\R^2)}^2
				\leq C_E^2M_{\overlap}C_{\hat{T}_0}^2h^{2\min\{s,p\}}\|u\|_{\HSo^s(\Omega)}^2\,,
	\end{aligned}
	\label{eq:uglobalhestimate}
\end{equation}
where $h=\sup_{K\in\mesh}h_K$ and $C_E=C_E(s,\Omega)$ is not dependent on $p$.
The same estimate holds for the variable $\bsq\in(\HSo^s(\Omega))^2$, and similar bounds (even without using extension operators and the finite overlap condition) hold for $\hat{u}\in\HSo_0^{\onehalf+s}(\bdry\mesh)$ and $\hat{q}_{\nml}\in\HSo^{-\onehalf+s}(\bdry\mesh)$.
Then, assume $M_F$ is independent of the family of meshes being considered, and choose the interpolant in \eqref{eq:convergencewithFortin} along with the estimates of the type in \eqref{eq:uglobalhestimate}, so that 
\begin{equation}
	\|\fku-\fku_h\|_\scU\leq\frac{\|b\|\,M_F}{\gamma}\|\fku-\Pi_{\scU^s}\fku\|_\scU
		\leq Ch^{\min\{s,p\}}\|\fku\|_{\scU^s}\,,
\end{equation}
where $C=C(p,s,\Omega)>0$, but is independent of the meshes being considered.
Moreover, if $M_F$ is $p$-independent, then in the $p$-asymptotic limit, the following $hp$-convergence estimate holds (see \eqref{eq:pdependenceofC}),
\begin{equation}
	\|\fku-\fku_h\|_\scU\leq\tilde{C}(\ln p)^2\frac{h^{\min\{s,p\}}}{p^s}\|\fku\|_{\scU^s}\,,
	\label{eq:hpconvergencefull}
\end{equation}
where $\tilde{C}=\tilde{C}(s,\Omega)$ is independent of $p$.
This concludes the results summarized in Theorem \ref{thm:convergence}.

\begin{remark}
Starting directly from the quasi-optimal error estimate in \eqref{eq:convergencewithFortin}, and avoiding interpolation, it is possible to get a better estimate for the variable $\hat{u}\in\HSo_0^{\onehalf+s}(\bdry\mesh)$ by using the results in \cite{babuskasuriapprox}, provided the triangulations $\scT(K)$ are quasi-uniform across all $K\in\mesh$ and all meshes being considered.
In that case, the $hp$-convergence estimate in \eqref{eq:hpconvergencefull} will have a $\ln p$ instead of $(\ln p)^2$.
\end{remark}
